# Discrete geometry and topology of entanglement of straight lines in 3-space

*Peter V Pikhitsa[1] and Stanislaw Pikhitsa[2]*

[1]Seoul National University, Seoul, Korea, peterpikhitsa@gmail.com; peter@snu.ac.kr

[2]Odessa National University, Odessa, Ukraine, aizetone@gmail.com

1 May 2020

**Abstract**

We propose an unexpected twist to description of the geometry and topology of configurations of $n$ straight lines considered as a whole 3D entity (because the lines are inseparably linked pairwise while having linking numbers ½ or -½) and named $n$-cross. Our theory stems from our work on configurations of mutually touching straight cylinders but, along with the previously introduced Ring matrix (that controls the encaging of each line by other lines), we now introduce fundamental *direction* 3D matrices (whose entries 0, 1, and -1 are signs of mixed products of line orientation vector triples). Discrete motion/connection combination principle established in the space of Ring and direction matrices (forming a groupoid and resembling moves in Loyd's 15-puzzle game or Khovanov homology) allows one to discern topologically different configurations of lines with elementary methods and *without* link diagrams of knot theory. However, with the help of so-called *projection* 3D matrix we also integrated our matrix approach into the knot theory and established topological invariants for line entanglement in both approaches thus connecting 2D projections with 3D configurations. With Jones polynomials we show that an $n$-cross is a link of pairwise connected $n$ unknots in a topological sense. The known results of the knot theory for rigid isotopy of 6 and 7 lines are reproduced and a novel result for 8 lines is given. With our approach we reach nuances of the geometry of lines never investigated before. It may find applications in Algebra, Discrete Geometry and Topology, and Quantum Physics.

*"A hidden connection in stronger than an obvious one"*, Heraclitus

**Introduction**

Some time ago we developed a manifestly 3D approach to solve the problem of finding configurations of mutually touching infinite straight cylinders in 3D (see [1,2] and references therein), which are straight "thick lines". Our classification of configurations (called $n$-knots when all the cylinders are in mutually touching, and $n$-crosses, with arbitrary positions of cylinders) was based on two matrices: a chirality matrix $P$ and a Ring matrix $R$ [1,2]. We developed a set of invariants that distinguished $n$-knots from $n$-crosses and found out that only one topologically unique 7*-knot, first discovered in [3], exists for cylinders of equal radii (along with its mirror image) [2]. No $n$-knots with equal radii of the cylinders can exist for $n > 7$. Yet for cylinders of arbitrary cross-section, $n$-knots are possible for $n < 11$ [1]. For example, a 9-knot made of equal cylinders with elliptical cross-sections is possible, as well as 10-knots with unequal elliptical cross-sections. Without the highly restrictive conditions of mutually touching, $n$-crosses exist for any $n$, but they become non-trivially topologically entangled only for $n > 5$ [4]. This entanglement might be fundamental and common in Nature, and calls for simple methods to control it.

Here we show that we can apply $P$, $R$ and some other 3D matrices (matrix-valued vectors, like the direction matrix $\widehat{N}$ introduced below) to solve the problem of entanglement of straight lines in 3D [4] in a discrete way. To provide a discrete analogy with Witten's theory of Jones polynomials in 3D, one can say that $P$ plays the role of the term of the product of Wilson loops in the integral of the Chern-Simon action over the gauge fields, and $\widehat{N}$ plays the role of the Chern-Simon action part. Our approach is based on matrix algebra, it is essentially discrete and is different from the knot theory because it can do without link diagrams, though it may use some of knot theory results. We investigate the topology of the entangled configuration exploring a discrete connection rule in the space of matrices that represent the finite configuration space of $n$-crosses. The latter seems to resemble the Khovanov homology approach. We found an intricate relation between a 3D $n$-cross and its 2D projection, thus viewing Witten's solution of the Atiyah problem of 3D interpretation of the knot theory, in a discrete way. Earlier, the problem of line entanglement was solved with the Jones polynomial approach, modified for $RP^3$ space in [4] and based on 2D projections of lines. This purely topological approach has little contact with the geometry of straight lines save that the straight line is not homologous to zero, unlike the loops that are usually considered in the theory of links and knots.

A configuration of straight lines, an $n$-cross, as observed from afar, looks like a point source of lines issuing from it (Fig. 1a), quite similar to Faraday lines issuing from a point charge, but its "core" (magnified in Fig. 1b) reveals complex inner structure and thus possesses inner degrees of freedom. One can see that the core consists of lines that "miss" each other and may produce entanglement which cannot be disentangled without line crossing or without a line being exactly parallel to another one. Below we will show how to define the core geometry precisely.

(a) (b)

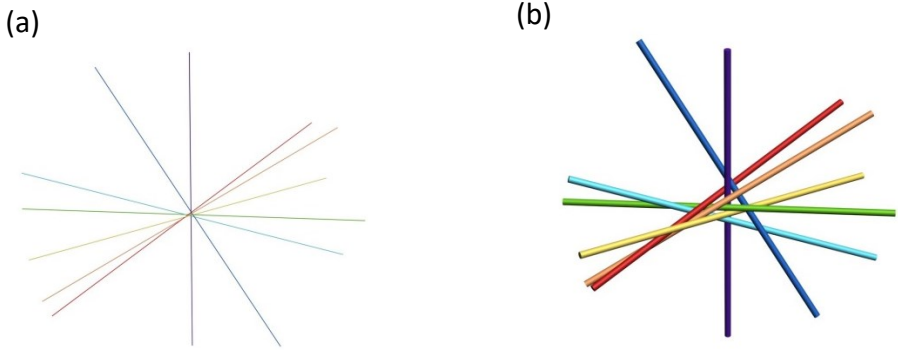

Fig. 1. An image of a 7*-knot (the configuration of 7 mutually touching equal cylinders first found in [3]) (a) from afar; (b) a close-up image of the core.

Unlike traditional knots and links, a configuration of oriented straight lines is inherently spin-like or "fermionic" in nature which is related to the $RP^3$ geometry of straight lines being not homologues to zero. Indeed, unlike a link of two closed loops where the Gauss linking number changes from 1 to 0 while one loop crosses the other, even a 2-cross is already non-trivial as far as its Gauss linking number is never 0 or 1. The oriented lines are described with $\boldsymbol{\gamma}_i(t_i) = \boldsymbol{n}_i t_i + \boldsymbol{v}_i$, where $\boldsymbol{n}_i$ is the unit vector of the line direction (note that in what follows we always enumerate lines starting from number $i = 0$) , and $\boldsymbol{v}_i = (x_i, y_i, 0)$ is a vector in the horizontal plane while the latter is punctured by the line at $t_i = 0$. Let us arrange the set $G$ of the line parameters of a configuration of an $n$-cross in the form of four vectors: $(G_0)_i = \theta_i$; $(G_1)_i = \varphi_i$; $(G_2)_i = x_i$; $(G_3)_i = y_i$ , so that with spherical angles $\theta_i$ , $\varphi_i$ the unit vector $\boldsymbol{n}_i = (\sin((G_0)_i)\cos((G_1)_i), \sin((G_0)_i)\sin((G_1)_i), \cos((G_0)_i))$ and $\boldsymbol{v}_i = ((G_2)_i, (G_3)_i, 0)$ .

The Gauss linking integral for two oriented straight lines $\boldsymbol{\gamma}_i, \boldsymbol{\gamma}_j$:

$$lk(\boldsymbol{\gamma}_i, \boldsymbol{\gamma}_j) = \frac{1}{4\pi} \int_{-\infty}^{\infty} \int_{-\infty}^{\infty} \frac{(\dot{\boldsymbol{\gamma}}_i(s), \dot{\boldsymbol{\gamma}}_j(t), \boldsymbol{\gamma}_i(s) - \boldsymbol{\gamma}_j(t))}{|\boldsymbol{\gamma}_i(s) - \boldsymbol{\gamma}_j(t)|^3} ds dt \qquad (1)$$

is either 1/2 or -1/2 [6] and changes to -1/2 or 1/2 , correspondingly, when one straight line crosses the other or when the configuration is mirrored (Fig. 2). One can say that an $n$-cross is always pairwise linked and should be considered (and described) as one whole entity. Below we will give a clear picture of it.

The sign of this link number becomes the element $P_{i,j}$ of the symmetric *chirality* matrix $P = \|P_{i,j}\|$ [1,2] with a zero diagonal and entries 1 and -1 that characterizes an $n$-cross:

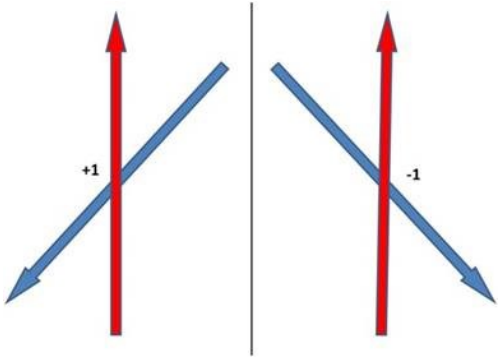

Fig. 2.

$$P_{i,j} = sign[lk(\pmb{\gamma}_i, \pmb{\gamma}_j)] \equiv sign[(\pmb{n}_i, \pmb{n}_j, \pmb{v}_i - \pmb{v}_j)]. \quad (2)$$

As an example, this matrix can be easily obtained directly from the 2D projection in Fig. 1b using the rules of Fig. 2 with a proper choice of the line directions to get, for example,

$$P = \begin{pmatrix} 0 & +1 & +1 & +1 & +1 & +1 & +1 \\ +1 & 0 & +1 & +1 & +1 & -1 & +1 \\ +1 & +1 & 0 & -1 & -1 & -1 & +1 \\ +1 & +1 & -1 & 0 & -1 & +1 & +1 \\ +1 & +1 & -1 & -1 & 0 & +1 & -1 \\ +1 & -1 & -1 & +1 & +1 & 0 & +1 \\ +1 & +1 & +1 & +1 & -1 & +1 & 0 \end{pmatrix} \quad (3)$$

so that its determinant is $\det(P) = -18$.

Another matrix which is needed to completely characterize an $n$-cross/knot is a Ring matrix $R$. It is defined as follows: non-diagonal entry $R_{i,j}$ is the number of triangles that encircle $i$-th line and contain $j$-th line as a side of the triangles. Its diagonal entries are zero. An example of $R$ of the 7*-knot from Fig. 1 is

$$R = \begin{pmatrix} 0 & 1 & 1 & 4 & 1 & 1 & 4 \\ 4 & 0 & 4 & 4 & 4 & 4 & 4 \\ 0 & 0 & 0 & 0 & 0 & 0 & 0 \\ 0 & 0 & 0 & 0 & 0 & 0 & 0 \\ 1 & 1 & 4 & 4 & 0 & 1 & 1 \\ 1 & 1 & 1 & 1 & 4 & 0 & 4 \\ 0 & 0 & 0 & 0 & 0 & 0 & 0 \end{pmatrix} . \quad (4)$$

Rows of zeroes indicate a "free" line that can be parallel translated to infinity without moving any other line. One can also build a Ring vector by summing up the numbers in each row and then divide each sum by 3. Its entries indicate the number of triangles that encircle a given line. We utilized some properties of the Ring matrix in previous papers on mutually touching cylinders and demonstrated its importance [1,2].

The Ring matrices have a far-reaching linear property that helps analyzing sub-configurations. For example, an $n$-cross has $n$ sub-configurations of $(n-1)$-crosses which Ring matrices are $R_{n-1}^{(i)}$, where $i$ indicates that $i$th line is omitted from the $n$-cross. If one adds to each of the Ring matrices of the latter sub-crosses corresponding row and column of zeroes and sums the matrices up, then one obtains the Ring matrix of an $n$-cross ($n > 4$):

$$R_n = \frac{1}{n-4} \sum_{i=0}^{n-1} R_{n-1}^{(i)}. \quad (5)$$

**Quantizing configurations and topology of n-crosses**

Now we show that there is a fundamental 3D *direction* matrix $\widehat{N}$ defined by the directions of the lines which is a "square root" of the Ring matrix in a way that there is a fundamental relation, valid for $n$-crosses:

$$R(P, \widehat{N})_{j,i} = \frac{1}{8}\left\{n(n-2) + 2P_{i,j}\left(\widehat{N}_j{}^2 P\right)_{i,j} - \left[\left(\widehat{N}_j P\right)_{i,j}\right]^2\right\}(1 - \delta_{i,j}), \quad (6)$$

where $\delta_{i,j}$ is the Kronecker delta,

$$\left(\widehat{N}_i\right)_{j,k} = sign(\boldsymbol{n}_i[\boldsymbol{n}_j \times \boldsymbol{n}_k]) \quad (7)$$

is the entry of the 3D direction matrix. We give the derivation of Eq. (6) in **Appendix 1**.

All $\widehat{N}_i$ have the same eigenvalues. For example, for $n = 6$ the characteristic equation for $\widehat{N}_i$ reads

$$\widehat{N}_i{}^2\left(\widehat{N}_i{}^4 + 10\widehat{N}_i{}^2 + 5\hat{I}\right) = 0,$$

where $\hat{I}$ is the identity matrix. Another distinguishing property of matrices $\widehat{N}_i$ is that each of them can be transformed by row/column permutations and row/column sign change into a unique form that can be called triangular. This form has all +1s above the zero diagonal and all -1s below (or *vice versa*), except the $i$th row/column which is filled with zeroes. This property reflects the possibility to arrange real directed lines in such a way that they look like a "fan" with arrows directed all in one half-plane being viewed along the $i$th line.

From Eq. (7) it is clear that the mixed product of vectors is anti-symmetric and produces one zero row/column in each square matrix component. A sign-switching for an entry happens for those triples where the co-planarity changes when lines move. For example, for a 6-cross, for a triple of lines $i = 0, j = 3, k = 4$ passing though co-planarity, three of the 6 matrices $\widehat{N}_i$, that is $\widehat{N}_0, \widehat{N}_3$, and $\widehat{N}_4$, change as

$$\widehat{N}_0 = \begin{pmatrix} 0 & 0 & 0 & 0 & 0 & 0 \\ 0 & 0 & -1 & +1 & -1 & +1 \\ 0 & +1 & 0 & -1 & 1 & -1 \\ 0 & -1 & 1 & 0 & \mathbf{1} & -1 \\ 0 & +1 & -1 & \mathbf{-1} & 0 & +1 \\ 0 & -1 & +1 & +1 & -1 & 0 \end{pmatrix} \rightarrow \widehat{N}_0' = \begin{pmatrix} 0 & 0 & 0 & 0 & 0 & 0 \\ 0 & 0 & -1 & +1 & -1 & +1 \\ 0 & +1 & 0 & -1 & 1 & -1 \\ 0 & -1 & 1 & 0 & \mathbf{-1} & -1 \\ 0 & +1 & -1 & \mathbf{1} & 0 & +1 \\ 0 & -1 & +1 & +1 & -1 & 0 \end{pmatrix}, \quad (8)$$

where prime indicates the switched direction matrix. As we said the direction matrix is not arbitrary in distribution of +1 and -1. It has definite signatures of belonging to real line configurations. For example, for a 6-cross it has four classes which can be defined as a function:

$$class(\widehat{N}) = tr\left(\frac{1}{\sum_{i=0}^{n-1} \widehat{N}_i^{\,2}}\right) \quad (9)$$

where possible zero eigenvalues of the matrix $x = \sum_{i=0}^{n-1} \widehat{N}_i^{\,2}$ in denominator are eliminated. Matrix $x$ satisfies equations

$$x(x+40)(x+4)^2(x+36)^2 = 0;$$

$$(x+4)(x+36)(x^2+40x+80)^2 = 0;$$

$$(x^2+40x+80)^2(x^2+40x+128) = 0;$$

$$(x^2+40x+80)^3 = 0,$$

so that Eq. (9) gives values for classes -0.580(5); -1.2(7); -1.3125; -1.5, correspondingly. Note that for $n < 6$ there is only 1 class for all configurations. For $n = 5$ this class equals -2.05.

We found that a direction matrix $\widehat{N}$ from a real $n$-cross configuration satisfies an identity looking quite similar to Eq. (7)

$$(\widehat{N}_i)_{j,k} = sign\{tr(\widehat{N}_i[\widehat{N}_j, \widehat{N}_k])\}, \quad (7a)$$

where the square brackets mean a commutator.

The switching event produces a rigid isotopy change in an $n$-cross with the corresponding change in the Ring matrix according to Eq. (6). The most important is that such a change gives us a possibility to establish a *combination principle* or a "Golden rule" for a *connection* which defines a correct switching/morphism between adjacent configurations, that reflects the continuity of motion in space and *discrete/quantum* topology by following the difference in the Ring matrices before and after switching, e.g:

$$R(P,\widehat{N}) - R(P,\widehat{N}') =$$

$$\begin{pmatrix} 0 & 0 & 0 & 0 & 0 & 0 \\ 3 & 0 & 1 & 3 & 1 & 1 \\ 1 & 3 & 0 & 1 & 3 & 1 \\ 0 & 0 & 0 & 0 & 0 & 0 \\ 0 & 0 & 0 & 0 & 0 & 0 \\ 2 & 2 & 4 & 2 & 2 & 0 \end{pmatrix} - \begin{pmatrix} 0 & 0 & 0 & 0 & 0 & 0 \\ 3 & 0 & 1 & 3 & 1 & 1 \\ 1 & 3 & 0 & 1 & 3 & 1 \\ 3 & 1 & 1 & 0 & 3 & 1 \\ 0 & 0 & 0 & 0 & 0 & 0 \\ 2 & 2 & 4 & 2 & 2 & 0 \end{pmatrix} = \begin{pmatrix} 0 & 0 & 0 & 0 & 0 & 0 \\ 0 & 0 & 0 & 0 & 0 & 0 \\ 0 & 0 & 0 & 0 & 0 & 0 \\ -3 & -1 & -1 & 0 & -3 & -1 \\ 0 & 0 & 0 & 0 & 0 & 0 \\ 0 & 0 & 0 & 0 & 0 & 0 \end{pmatrix}. \quad (10)$$

Note, that in this example only the line 3, being sandwiched between lines 0 and 4 at the moment of the co-planarity switching, produces the single row in the Ring matrix difference of Eq. (10). We use this property to characterize the "connected cluster" (a groupoid) of rigid isotopy configurations for any given $P$ that are connected in the sense of Eq. (10). The configurations that are not connected cannot show one row at a single switching of the direction matrix $\widehat{N}$. This holds for any $n$-crosses and completely defines the rigid isotopy of any configuration of lines.

In the example of Eq. (10) the only non-zero row consists of (1 1 1 3 3 0) [along with its negative (-1 -1 -1 -3 -3 0)], yet, for a 6-cross there is another row (1 1 1 1 -1 0). Generally, for an $n$-cross the number of rows of different content types is $\left[\frac{n-3}{2}\right] + 1$ where the square brackets mean the integer part. For a 4-cross the single row is (1 1 1 0), for a 5-cross the 2 rows are (1 1 2 2 0) and (1 -1 0 0 0); for a 7-cross the 3 rows are (1 1 1 1 4 4 0), (1 1 1 -1 2 2 0), and (1 1 -1 -1 0 0 0); for an 8-cross the 3 rows are (1 1 1 1 1 5 5 0), (1 1 1 1 -1 3 3 0), and (1 1 1 1 1 -1 -1 0), and so on. Here numbers different from 1, stand in places with the numerals of the two lines in the sandwich of the switching triple. For example, in case of Eq. (10) the places are the 0[th] and 4[th] in the row, where -3s stand.

Let us make some remarks. Our approach, by establishing natural connection rules (forming a groupoid and resembling Khovanov homology) in the discrete space of configurations, gives a new twist to the straight line entanglement problem providing a groupoid calculus of evaluation of the entanglement. It is remarkable that Heisenberg discovered quantum mechanics by considering a groupoid of transitions for the hydrogen spectrum (to which Schwinger gave an algebra), rather than the usually considered group of symmetry of an individual state [9]. Here the Ring matrix plays a role of a Hamiltonian which registers the changes in states (presented by the configurations of lines) and determines the adjacent ones. Its additive nature expressed in Eq. (5) also supports this conclusion. Being explicitly 3D, our approach is different from the methods of the knot theory with its link diagrams, where the connection between changing configurations is established with the Redeimeister moves applied to the projection and the skein relations based on them, leading to Jones polynomials. We feel that our approach can be related to Vassilyev invariants (not only through Jones polynomials), because line intersections just produce zeroes in corresponding entries of the chirality matrix which creates a rich set of additional invariants. However, this is out of scope of the current work.

As an example, we applied the fundamental rule of Eq. (10) of rigid isotopy moves directly in 3D to investigate the entanglement of 6-crosses, the first non-trivial case of straight line entanglement studied in [4-7]. To control the configurations in connected clusters we introduce a configuration invariant

$$Inv(P, \widehat{N}) = \text{tr}\left(\frac{1}{\sum_{i=0}^{n-1} \widehat{N}_i^2 - \frac{P}{2}}\right). \quad (11)$$

The invariant of Eq. (11) based on the direction matrix distinguishes geometrically different configurations (compare with previously introduced invariants based on the Ring matrix [2]).

Some chirality matrices $P$ while having the same determinant can be distinguished by their set of eigenvalues. However, sometimes it is not enough, because in general matrices might be not similar even having the same set of eigenvalues. Here we introduce a 3D matrix with entries

$$[T3(P)_i]_{j,k} = P_{i,j} P_{j,k} P_{k,i} \quad (12)$$

and use, for example, a number $InvP(P) = tr(\sum_i (T3(P)_i + \hat{I}/2)^{-1})$ to characterize the chirality matrices. With a given determinant of the chirality matrix for a 6-cross $|P| = 11,19,-21,27,-29,-13,-45,-5,-5^*,-125$ (and their mirror configurations marked with letter m in Table 1) we explored the connected clusters of all possible discrete configurations in clusters. The results are given in Table 1 along with $InvP(P)$ and the values of Jones polynomials $J_D(|P|, i, a)$ (see below) calculated at $a = 0.8$ ($i$ just enumerates topologically different clusters of the same $|P|$). The exceptional three configurations with $|P| = -45, -29, -13$ are equivalent in rigid isotopy to their mirror configurations and are called specular in [7]. The total number of isotopy different configurations is 19 as was proved in [4,6,7,8]. The sum of numbers in the 4th column of all possible discrete configurations of 6 straight lines in 3D space is 11618.

|   | $|P|$ | $InvP(P)$ | $DJ(|P|, i, 0.8)$ | Cluster size | $\sum[Inv(P, \widehat{N})]$ |
|---|---|---|---|---|---|
| 1 | -125 | 21.47368 | 81.95805 | 112 | -160.15626 |
| 2 | -125** | 21.47368 | 82.84623 | 112 | -161.01855 |
| 3 | -45 | 8.36522 | 85.95424 | 2256 | -3556.33638 |
| 4 | -29 | 0.93146 | 89.17975 | 1835 | -2842.88768 |
| 5 | -21 | 1.90476 | 99.95387 | 448 | -836.91785 |
| 6 | -21m | -1.80088 | 82.28432 | 448 | -668.07359 |
| 7 | -13 | -0.82759 | 91.67157 | 187 | -301.51497 |
| 8 | -5 | 14.46046 | 94.914 | 2100 | 142.24454 |

| 9  | -5m    | 0.51093   | 80.21999  | 2100 | 2266.87334 |
| 10 | -5*    | -34.66667 | 64.04794  | 16   | -7.52044   |
| 11 | -5*m   | 26.28571  | 162.15833 | 16   | -40.17786  |
| 12 | 11     | -17.99234 | 72.137    | 161  | -224.09236 |
| 13 | 11m    | 14.09524  | 125.74818 | 161  | -250.22019 |
| 14 | 19     | -14.13903 | 76.55007  | 635  | -981.15722 |
| 15 | 19m    | 14.00147  | 108.91363 | 635  | -1426.52974 |
| 16 | 27     | -10.28571 | 81.85074  | 149  | -128.16471 |
| 17 | 27m    | 13.90769  | 94.24605  | 149  | -219.1273  |
| 18 | 27**   | -10.28571 | 83.23852  | 49   | -62.15941  |
| 19 | 27**m  | 13.90769  | 93.67761  | 49   | -455.79888 |

Table 1.

The 5th column gives a sum of all invariants $\sum[Inv(P,\widehat{N})]$ for each cluster, thus providing a topological invariant for rigid isotopy. Then for any given configuration, described by $P,\widehat{N}$ one can find to which cluster it belongs by reconstructing all the cluster invariants using the "one-row" rule of Eq. (10). After exhausting all the possibilities for obtaining any new $Inv(P,\widehat{N})$ on the way, one gets essentially 3D topological invariant as the sum of all of the different ones:

$$\mathfrak{G} = \sum[Inv(P,\widehat{N})].$$

Non-trivially entangled configurations are those with $|P| = 27$. There are two topologically different clusters of topologically connected configurations: one with 149 configurations (and its mirror marked 27m) and the other one with 49 configurations (marked 27** and its mirror 27**m). We put into Appendix 2 Table A2.1 of the invariants $Inv(P,\widehat{N})$ for all 49 configurations of cluster 27**. The Table also shows the adjacent configurations in the cluster so that the topology of the cluster can be additionally characterized by its Euler characteristic. The mirror cluster 27**m has different set of values of $Inv(P,\widehat{N})$. Further we will calculate a Jones polynomial for the lines (first established in [5]) along with a novel one that we call $J_M$ – polynomial; comparison between two polynomials reveals the difference in entanglement for both clusters.

As a corollary, we directly, by inspecting all 224 configurations with $|P| = -125$, proved our previous conjecture [1] that only two of the configurations (with mirrored ones) can allow all mutually touching cylinders and thus can be a 6-knot when $|P| = -125$. The uniqueness of these configurations was a cornerstone in our proof that there is a bottleneck preventing mutual touching of more than 10 arbitrary cylinders in 3D [1].

**Plane projection of $n$-crosses**

Let us describe relations between the plane projections of line configurations with 3D structures. In fact this is a discrete version of Witten's idea to connect knots with Jones polynomials built on 2D link diagrams. Projections give a possibility to use heavy artillery of the knot theory. The groupoid rules of Eq. (10) give us an alternative method to describe the topology. We also manage to characterize the geometry of $n$-crosses to such an extent which is not possible for any knot theory methods.

Introduce a matrix valued vector $prM(G, \boldsymbol{U})_i$ which we call a projection matrix similar to direction matrix $\widehat{N}_i$, which is now based on a 2D projection of the straight line configuration $G$ along some vector $\boldsymbol{U}$ onto a plane. Its entry $[prM(G, \boldsymbol{U})_i]_{j,k}$ is defined as shown in Fig. 3 with a simple rule: the value is +1 if the direction from the crossing point $(i, j)$ to the crossing point $(i, k)$ coincides with the direction of line $i$ and is -1 when opposite as it is in Fig. 3. Therefore $[prM(G, \boldsymbol{U})_i]_{j,k} = -1$ here. Yet, $[prM(G, \boldsymbol{U})_j]_{k,i} = 1$ after index permutations which one would not expect from $[\widehat{N}_i]_{j,k}$: the latter stays the same for any cyclic index permutation.

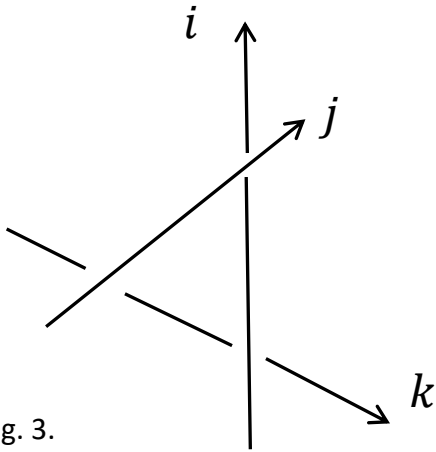

Fig. 3.

This deceivingly simple definition produces a solution of the problem how to connect 3D configurations with their 2D projections.

There are general properties of $prM(G, \boldsymbol{U})$. First,

$$R(UU, prM(G, \boldsymbol{U})) \equiv \widehat{0}, \quad (13)$$

(this is the general property of a 3D matrix which vector components are similar to anti symmetrized matrix $UU$, that is triangular) where the square matrix $UU$ has a zero diagonal and all other entries are 1. We will use Eq. (13) to define the geometry of the inner domain of $n$-cross later on. Second, by observing all possible projections of three lines one can obtain a general relation for $prM(G, \boldsymbol{U})$:

$$[prM(G,\boldsymbol{U})_i]_{j,k}[prM(G,\boldsymbol{U})_j]_{k,i} = -O_{i,k}O_{j,k}P_{i,k}P_{j,k} \, , \quad (14)$$

where $i \neq j$ and we introduced an anti-symmetric overlapping matrix $O_{i,k}$ which entry is +1 when the projection line $i$ overpasses the line $k$ and -1 otherwise. For example, from Fig. 3 one can obtain $O_{i,k} = +1$, $O_{j,k} = +1$, $P_{i,k} = -1$, and $P_{j,k} = -1$, while $[prM(G,\boldsymbol{U})_i]_{j,k} = -1$ and $[prM(G,\boldsymbol{U})_j]_{k,i} = 1$ so that Eq. (14) is correct. If we multiply Eq. (14) by $[prM(G,\boldsymbol{U})_k]_{i,j}$ we will obtain a cyclic-invariant entry $[\widehat{N}_i^{\,c}]_{j,k}$:

$$[\widehat{N}_i^{\,c}]_{j,k} = [prM(G,\boldsymbol{U})_i]_{j,k}[prM(G,\boldsymbol{U})_j]_{k,i}[prM(G,\boldsymbol{U})_k]_{i,j} \, . \quad (15)$$

One can apply Eq. (15), directly obtaining $[\widehat{N}_i^{\,c}]_{j,k}$ from the projection in Fig. 3. Indeed, the triple of line indexes $i,j,k$ outside the triangle naturally determines a vortex direction (clockwise in Fig. 3). Then this direction is opposite to the direction of line $i$, agrees with line $j$, and disagrees with line $k$ to give $(-1)(1)(-1) = 1$ which coincides with $[\widehat{N}_i^{\,c}]_{j,k} = 1$ from Eq. (15).

Because of the importance of $\widehat{N}^c$, which connects 2D projections and 3D configurations as we show below, we introduce a function that expresses Eq. (15) in a short form:

$$\widehat{N}^c = D3(prM(G,\boldsymbol{U})) \, . \quad (16)$$

Alternatively,

$$[\widehat{N}_i^{\,c}]_{j,k} \equiv [\widehat{N}_k^{\,c}]_{i,j} = -O_{i,k}O_{j,k}P_{i,k}P_{j,k}[prM(G,\boldsymbol{U})_k]_{i,j} \, . \quad (17)$$

Let us introduce a function

$$[D2(A,B)_k]_{i,j} = A_{i,k}A_{j,k}[B_k]_{i,j} \, , \quad (18)$$

so that Eq. (17) can be rewritten as

$$prM(G,\boldsymbol{U}) = -D2(P \otimes O, \widehat{N}^c) \, , \quad (19)$$

where the circled times sign means the direct product. As an exercise one can show that there is an identity

$$R(P, \widehat{N}) \equiv R(UU, D2(P, \widehat{N})) \, .$$

With the help of Eq. (19) we can rewrite Eq. (13) as

$$R\big(UU, D2(P \otimes O(G,\boldsymbol{U}), \widehat{N}^c)\big) = R(UU, D2(P \otimes O(G,\boldsymbol{U}), D3(prM(G,\boldsymbol{U})))) \equiv \widehat{0}. \quad (20)$$

From Eqs. (15) and (17) one can derive the chirality matrix $P$ (up to a sign) as a function of $prM(G, \boldsymbol{U})$ and $O(G, \boldsymbol{U})$ (or *vice versa*). For example, Eq. (17) can be rewritten as

$$(P \otimes O)^2 = \sum_k \widehat{N}_k^{\,c} \otimes prM(G, \boldsymbol{U})_k - (n-1)\hat{I}, \quad (21)$$

where $\widehat{N}_k^{\,c}$ is taken from Eq. (15). While extracting the square root of a matrix in Eq. (21) is a cumbersome procedure, we found a function that directly calculates $H = P \otimes O$ through $prM(G, \boldsymbol{U})$ (see **Appendix 5**).

Another discrete connection of the 2D projection of the initial 3D configuration with the direction matrix $\widehat{N}$ comes from Eq. (20), where one can replace $\widehat{N}^c$ by $\widehat{N}$ and the equation still holds for any arbitrary projection vector $\boldsymbol{U}$

$$R\big(UU, D2(P \otimes O(G, \boldsymbol{U}), \widehat{N}(G))\big) = \hat{0} \; . \quad (22)$$

This indicates that the original $\widehat{N}(G)$ is the "kernel" of projections.

It happens that for projections matrix $\widehat{N}^c$ belongs to the same set of classes of Eq.(9) as the direction matrix $\widehat{N}(G)$ does. Moreover, there is the most remarkable property of $\widehat{N}^c$:

*If one calculates $Inv(P, \widehat{N}^c)$ then one gets a value lying within the full set of values for a connected cluster (including $Inv(P, \widehat{N}(G))$) of the initial 3D configuration $G$.*

That means that by rotating the projection vector $\boldsymbol{U}$ we obtain exclusively the values of invariants that belong to a connected cluster of an $n$-cross and thus keep the information of the topology of the cluster. However, not all the cluster is covered with these values because of an excess connectivity demanded by the conditions of the projection. One may say that the projections define a sub-groupoid, unlike the entire groupoid that can be found with spanning the cluster with connection rules of Eq. (10).

Also there is a caveat when dealing with projections. For an exceptional chirality matrix with $|P| = -125$ for a 6-cross, there are two clusters as shown in Table 1 and $prM(G, \boldsymbol{U})$ gives $Inv(P, D3(prM(G, \boldsymbol{U})))$ which is always the invariant of the other cluster, different from the cluster to which $Inv(P, \widehat{N}(G))$ of the initial 3D configuration belongs. Still the clusters remain quite separable and we can distinguish them from each other. As a sub-structure, this exceptional 6-cross with $|P| = -125$ appears in many $n$-crosses ($n > 6$), yet it does not violate the possibility to distinguish topologically different clusters.

Physically it is clear why the topology is preserved. The projection rotation does not change entanglement of the straight lines; one can say that it realizes Reidemeister move III for the lines when the projection image changes as the entire $n$-cross rotates. For a continuous change in $U$ sweeping the sphere, switching events in the direction matrix happen automatically, and the Ring matrix changes according to the one-row connection rule of Eq. (10) as well. All the configurations turn out to be connected in the discrete topology naturally. The domains/patches of equal invariants $Inv(P, \widehat{N}^c)$ on the sphere swept by $U$ provide a tessellation on the spherical surface.

Now we will give a concrete example how the projection works with the identification of the line configuration of a 6-cross provided in [4] (Fig. 4a) where we equipped the lines with numbers and directions.

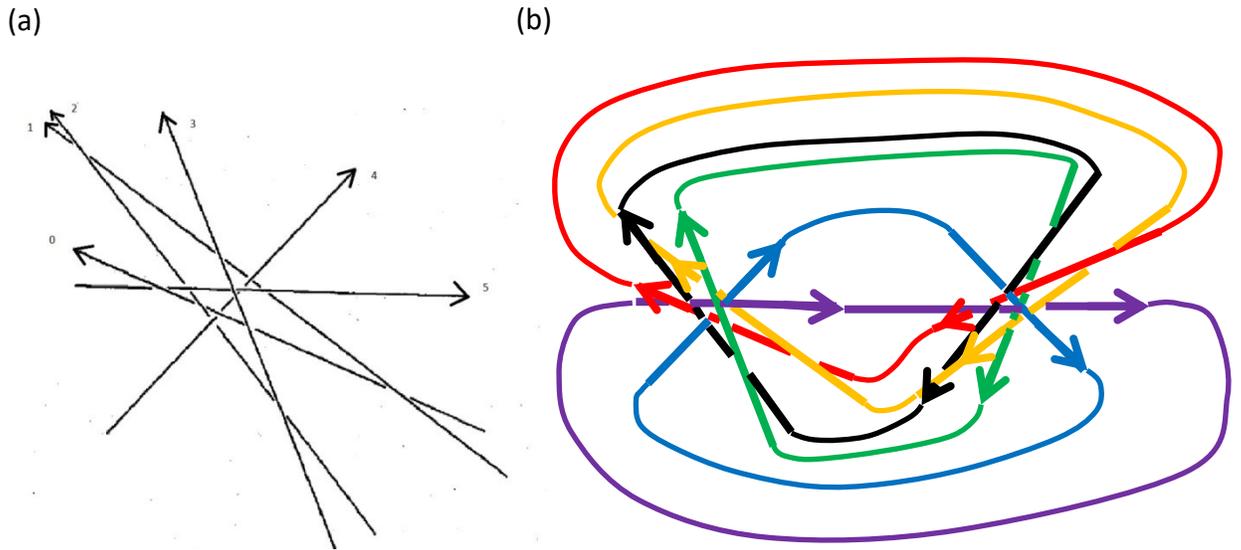

(a)    (b)

Fig. 4.

According to rules in Fig. 2 one can get the chirality matrix from Fig. 4

$$P27^{**} = \begin{pmatrix} 0 & +1 & +1 & +1 & -1 & -1 \\ +1 & 0 & +1 & +1 & +1 & +1 \\ +1 & +1 & 0 & +1 & +1 & +1 \\ +1 & +1 & +1 & 0 & -1 & -1 \\ -1 & +1 & +1 & -1 & 0 & +1 \\ -1 & +1 & +1 & -1 & +1 & 0 \end{pmatrix} \quad (23)$$

which determinant is 27 and $InvP(P) = -5.82857$. Then one can get components of the projection matrix from Fig. 4

$$prM_0 = \begin{pmatrix} 0 & 0 & 0 & 0 & 0 & 0 \\ 0 & 0 & -1 & -1 & -1 & -1 \\ 0 & +1 & 0 & +1 & +1 & -1 \\ 0 & +1 & -1 & 0 & -1 & -1 \\ 0 & +1 & -1 & +1 & 0 & -1 \\ 0 & +1 & +1 & +1 & +1 & 0 \end{pmatrix}; prM_1 = \begin{pmatrix} 0 & 0 & -1 & -1 & -1 & -1 \\ 0 & 0 & 0 & 0 & 0 & 0 \\ +1 & 0 & 0 & +1 & +1 & +1 \\ +1 & 0 & -1 & 0 & +1 & +1 \\ +1 & 0 & -1 & -1 & 0 & +1 \\ +1 & 0 & -1 & -1 & -1 & 0 \end{pmatrix}; prM_2 = \begin{pmatrix} 0 & -1 & 0 & +1 & +1 & -1 \\ +1 & 0 & 0 & +1 & +1 & +1 \\ 0 & 0 & 0 & 0 & 0 & 0 \\ -1 & -1 & 0 & 0 & -1 & -1 \\ -1 & -1 & 0 & +1 & 0 & -1 \\ +1 & -1 & 0 & +1 & +1 & 0 \end{pmatrix};$$

$$prM_3 = \begin{pmatrix} 0 & -1 & +1 & 0 & -1 & -1 \\ +1 & 0 & +1 & 0 & +1 & +1 \\ -1 & -1 & 0 & 0 & -1 & -1 \\ 0 & 0 & 0 & 0 & 0 & 0 \\ +1 & -1 & +1 & 0 & 0 & -1 \\ +1 & -1 & +1 & 0 & +1 & 0 \end{pmatrix}; prM_4 = \begin{pmatrix} 0 & -1 & +1 & -1 & 0 & -1 \\ +1 & 0 & +1 & +1 & 0 & +1 \\ -1 & -1 & 0 & -1 & 0 & -1 \\ +1 & -1 & +1 & 0 & 0 & -1 \\ 0 & 0 & 0 & 0 & 0 & 0 \\ +1 & -1 & +1 & +1 & 0 & 0 \end{pmatrix}; prM_5 = \begin{pmatrix} 0 & -1 & -1 & -1 & -1 & 0 \\ +1 & 0 & +1 & +1 & +1 & 0 \\ +1 & -1 & 0 & -1 & -1 & 0 \\ +1 & -1 & +1 & 0 & -1 & 0 \\ +1 & -1 & +1 & +1 & 0 & 0 \\ 0 & 0 & 0 & 0 & 0 & 0 \end{pmatrix}. \quad (24)$$

Now from Eqs. (11),(16),(23) and (24) we get $Inv(P27^{**}, D3(prM)) = -1.740388404194$ which value is 21st in the Table A2.1 of **Appendix 2**, marked red. Again, if we started with the current $P27^{**}$ and $\widehat{N}^c = D3(prM)$ we could obtain the whole cluster with 49 elements by the connection rule of Eq. (10). Thus the projection gives us the complete information of the entanglement of the lines, prescribing to the configuration in Fig. 4 the topological invariant -62.15941 of the fifth column in Table 1.

As we said above, if, for a given $n$-cross configuration $G$, one rotates $U$ in all directions to cover the solid angle $4\pi$, still the projections cannot "cover" the cluster completely, which means that the number of invariants $Inv(P(G), D3(prM(G, U)))$ is always less than the total size of the cluster connected with the one-row rule of Eq. (10). This happens because the projections only imitate the real switching in the direction matrix, thus selecting only specific connections between the configurations. Yet, they can be used to characterize the geometry of an $n$-cross by extending the notion of the projection matrix as given below.

**Inside and outside of an n-cross**

One can generalize the projections to investigate the space volume between the lines. These "3D point projections" can be defined as follows. Let now $U$ be not a projection direction but a 3D coordinate vector of a point that itself issues projection rays. Then the "shadows" of lines $j$ and $k$ intersect with the shadow of line $i$ as in Fig. 3, so that we obtain a matrix-valued vector with entries $[prM3D(G, U)_i]_{j,k}$ defined as before for $[prM(G, U)_i]_{j,k}$. However, the class of $D3[prM3D(G, U)]$ may not be correct (not coinciding with a class of a real configuration) because the projecting point $U$ may be sandwiched in between the lines $j$ and $k$. In between means that $U$ lies between the planes, one of which contains line $j$ and is parallel to line $k$ and the other one contains line $k$ and is parallel to the line $j$, so that projections are in opposite directions. This can be corrected with a symmetric matrix $P3D(G, U)_{j,k}$ of zero diagonal which is -1 for the in-between case and 1 otherwise. The corrected projection matrix with entries

$$[prMN(G, \boldsymbol{U})_i]_{j,k} = [prM3D(G, \boldsymbol{U})_i]_{j,k}[T3(P3D(G, \boldsymbol{U}))_i]_{j,k} \quad (25)$$

(we used the function from Eq. (12)) produces $\widehat{N}^c = D3[prMN(G, \boldsymbol{U})]$ which has correct classes. We get more invariants $Inv\big(P(G), D3(prMN(G, \boldsymbol{U}))\big)$ of the cluster when $\boldsymbol{U}$ runs the whole bulk space but additional invariants come only from *internal* part of the n-cross (still in not enough quantity to cover the whole cluster). Indeed, when $\boldsymbol{U}$ is large enough (the projection point is far from the core of configuration), then we return to the case of plane projection as the projection rays become nearly parallel. One can provide a definition of the inside and outside domains of an $n$-cross by using the property of Eq. (13):

if $R(UU, prMN(G, \boldsymbol{U})) \equiv \hat{0}$  then $\boldsymbol{U}$ is outside; otherwise it is inside.    (26)

Alternatively, if $P3D(G, \boldsymbol{U}) - UU \equiv \hat{0}$  then $\boldsymbol{U}$ is outside; otherwise it is inside, that is sandwiched between at least two lines as described above.

**Jones topology invariants for $n$-crosses**

In parallel, we can apply the methods of the knot theory to confirm that the above described discrete topology complies well with the disentanglement procedure of skein relations that leads to two types of Jones polynomials which we call $J_D$-polynomials and $J_M$-polynomials. We use designation $J_D$ for Jones polynomials modified by Dobrotukhina for $RP^3$ in [4], and our novel polynomial $J_M$ created with closing lines into loops through "doubling" the configurations that can be recognized from the schematic in Fig.1 of [6].

 Let us define the skein relations for both polynomials and give elementary examples in Fig. 5 starting from 2-cross.

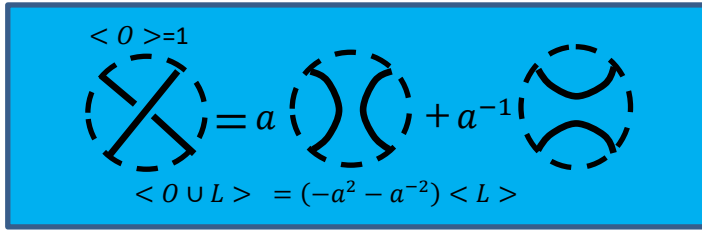
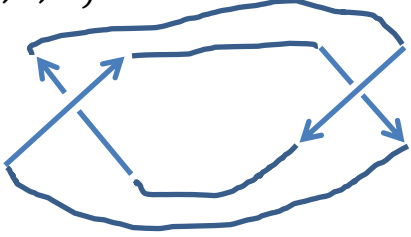
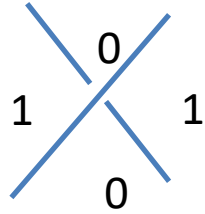

Fig. 5.

The skein rules are given in the rectangular in Fig. 5. On the left side in $J_M(1,0,a)$ the first argument means $|P| = 1$, the second one means the number of the cluster (we start numeration from 0).

The procedure of the doubling is the following. A copy of the 2-cross is flipped 180 degrees over the horizontal line, shifted to the right from the original and connected head-to-tail to the same arrows of the original. It produces the Hopf link $J_M(1,0,a) = -a^4 - a^{-4}$. The Hopf link reflects the nature of the pair of straight lines having the link number either 1/2 or -1/2 and never 0 so that they are always linked. Like the straight line, the Hopf link is not homologous to zero as well. Moreover, the structure is protected from passing the singularity when lines may become parallel. The Hopf link in Fig. 5 demonstrates this property. In fact, the $n$-link structure is in $S^3$ instead of initial $RP^3$ for lines. However, the doubling process applied to the diagram in Fig. 4a gives the diagram in Fig. 4b that exactly reproduces matrices from Eq. (24) and Eq. (25) when one uses the same rules for their calculations, only now applied to oriented circles instead of oriented straight lines.

Moreover, if one puts Fig. 4b on the surface of a large sphere near the North Pole and drags the right hand side doubled configuration to the South Pole then all the oriented circles will turn into large circles on the sphere but lifted a little off the sphere to make overpasses and underpasses with other circles. We illustrate it for a 3-cross in Fig. 6 where the part of the configuration (numbers with primes) on the right hand side from the dashed line was dragged rightwards until it comes to the South Pole and gives the corresponding three circles.

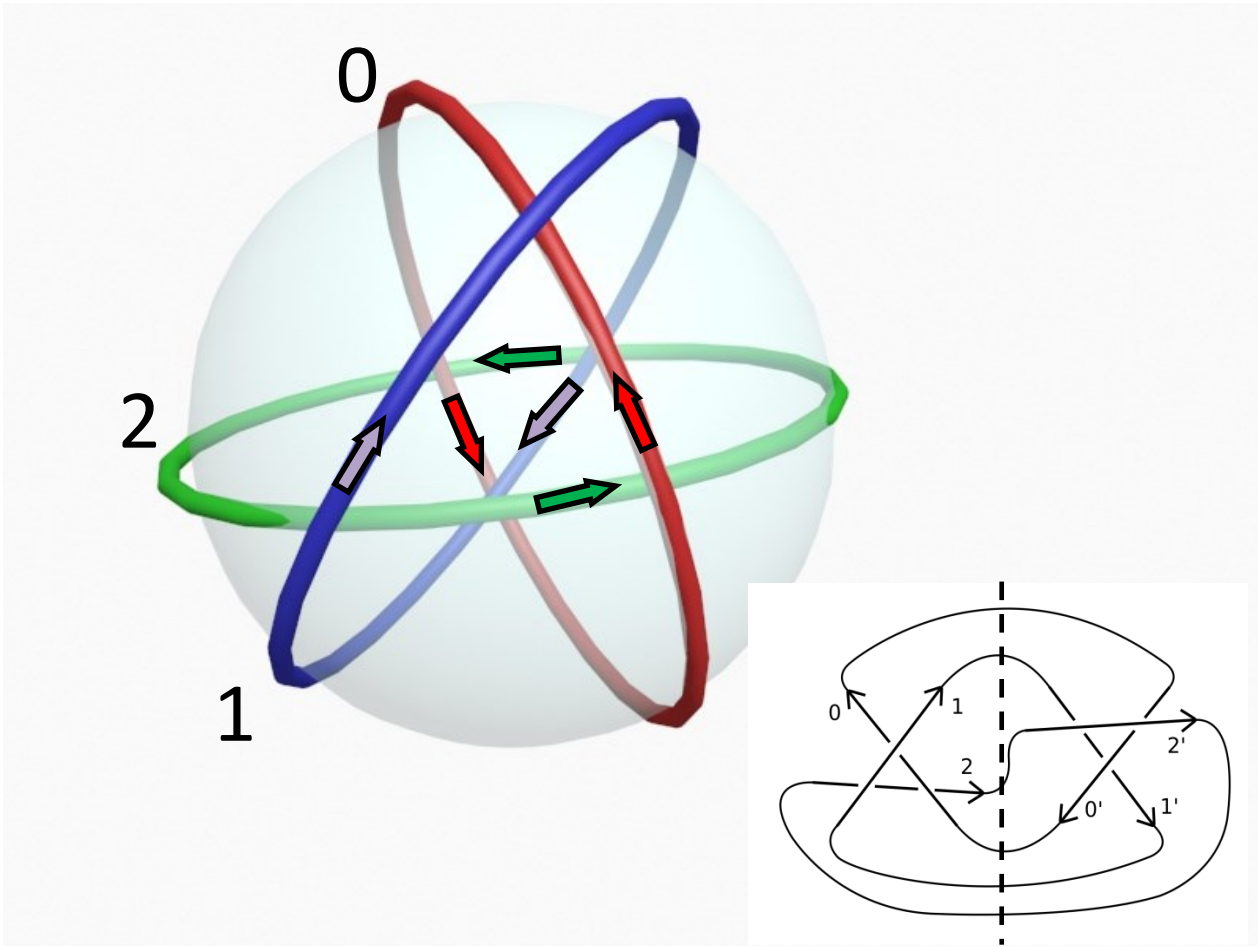

Fig. 6.

It is safe to say that *a configuration of n straight lines (an n-cross) is topologically equivalent to n rings (unknots) all linked pairwise*; it is a complete $n$-link (in analogy with a complete graph). It is clear that the famous Borromean rings are not in this set because they are not linked pairwise.

Note, that this equivalence to the circles explains why there exist connected clusters (groupoids) in $n$-crosses that are "impossible" to be realized with the straight line configurations but are quite legitimate in the domain of the complete $n$-link which is thus a generalization of an $n$-cross. Any further generalization would include unlinked circles which are formally described by a non-symmetric chirality matrix $P$ (to cover Borromean-type structures where this matrix is antisymmetric) which is an unexplored area for now. This type of arrangement of rings in 3D resembles the medieval armor called "Mail" and made of interlinked mesh of metal rings.

On the right side of Fig. 5 the identification of the opposite points in the projection of the 2-cross leads to identification of two domains 0 and 1 connected through infinity and thus protecting $RP^3$ geometry of lines. It is easy to see that the skein relations lead to two circles with factors $a$ and $a^{-1}$, so that $J_D(1,0,a) = a^1 + a^{-1}$. We give (for pedagogical reasons) detailed calculation of $J_M(2,0,a)$ and $J_D(2,0,a)$ for a 3-cross with $|P| = 2$ in **Appendix 3**.

Let us turn to a 6-cross. Both polynomials can be calculated with a designed computer program based on MathCad11. For the projection in Fig, 4 the result is:

$$J_D(27,0,a) = 22a + 15a^{-1} - a^3 - 12a^{-3} - 12a^5 - 13a^{-5} + a^7 + 10a^{-7} + 8a^9 + 15a^{-9} + 3a^{11} - 5a^{-13} + a^{-17} \quad (27)$$

which coincides with the result of [4,8] (in [8] it is given in Table 2 and labeled L among 19 rigid isotopy configurations) and gives $J_D(27,0,0.8) = 83.23852$ from Table. 1. The other polynomial reads

$$J_M(27,0,a) = 881 - 711a^4 - 963a^{-4} + 477a^8 + 913a^{-8} - 261a^{12} - 767a^{-12} + 97a^{16} + 541a^{-16} - 21a^{20} - 319a^{-20} - 3a^{24} + 141a^{-24} - 45a^{-28} + 9a^{-32} - a^{-36} . \quad (28)$$

For another (a topologically different cluster) 6-cross projection from [3] given in Fig. 7 we obtained

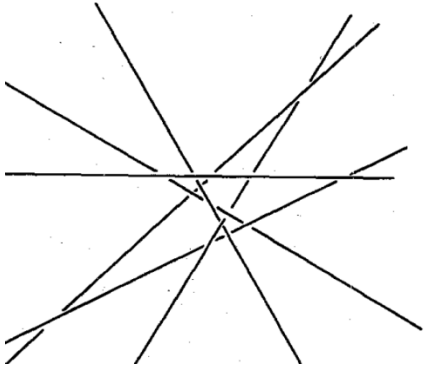

Fig. 7.

$$J_D(27,1,a) = 2a + 5a^3 + 3a^{-3} + 3a^5 + 7a^{-5} + 4a^{-7} + 2a^9 + 3a^{11} + a^{-11} + a^{13} + a^{-13} \quad (29)$$

which also coincides with the result of [4, 8] (in [8] this corresponds to the configuration in Table 3 labeled hc(125634)) and gives $J_D(27,1,0.8) = 81.85074$ from Table 1. The other polynomial reads

$$J_M(27,1,a) = -2 - 3a^4 - 3a^{-4} - 2a^8 - 4a^{-8} - 3a^{12} - 3a^{-12} - 2a^{16} - 4a^{-16} - a^{20} - a^{-20} - 2a^{24} - a^{28} - a^{-28} . \quad (30)$$

One may notice that Eqs. (29) and (30) are not independent! It is remarkable that they are connected with a simple equation

$-J_D(27,1,a^2)(-a^4 - a^{-4}) \equiv J_M(27,1,a)(-a^2 - a^{-2})$ (31)

where we deliberately left the signs as they are present in Hopf links and the skein rules in Fig. 5. For a 2-cross in Fig. 2 Eq. (31) is just an identity. It is easy to see using the polynomials for the 3-cross with $|P| = 2$ from **Appendix 3** that another relation holds:

$-J_D(2,1,a^2) \equiv J_M(27,1,a)$ .

For one of the 5-crosses with $|P| = 8$ still another relation occurs:

$-J_D(27,1,a^2) \equiv J_M(27,1,a)(-a^2 - a^{-2})$ .

Eq. (31), taken in a broader way as the equivalence

$-J_D(27,1,a^2) \sim J_M(27,1,a)$

(modulo Hopf link and unknot multipliers) works for any $n$-cross cluster which can be disentangled in a way that its connected groupoid contains a "trivial" configuration of lines with a zero Ring matrix, where each line is free to be translated to infinity. On the other side, Eqs. (27) and (28) do not satisfy $-J_D(27,0,a^2) \sim J_M(27,0,a)$ and thus present a non-trivially entangled $n$-cross cluster. Of all possible determinants for a 6 cross, Eq. (31) is fulfilled only for configurations with $|P| = 11, 19, -21, 27, -29, -13, -5^*$ (and their mirror ones) that can be disentangled.

We filled the third column of Table 1 with the values of $J_D(|P|, i, a)$ at $a = 0.8$, where, as we already mentioned, $i$ stays for a number of different clusters (for example, $i = 0$ and 1 for $|P| = 27$). One can notice that the invariants are in concordance with the last column which sums up all individual invariants in corresponding groupoids. Our Table 1 gives the same classification as Tables 1,2,3 from [8] but reflects a deeper view on the rigid isotopy of $n$-crosses because we also provide an essentially 3D invariant for rigid isotopy created outside knot theory.

In fact it is nearly impossible to span all possible configurations by generating lines in 3D randomly. We found that a discretization of the configurations makes it possible to use $J_D$ and $J_M$ polynomials to harvest rigid isotopy of all $n$-crosses with the help of the direction matrix $\widehat{N}$.

To do this let us take an $n$-cross with its direction matrix $\widehat{N}$. Now we have to produce a pseudo projection matrix $prM(\widehat{N})$ while not doing any real projection. In fact, we just simulate a would be projection in the vicinity of the direction of one straight line, while using the known entries of $\widehat{N}$ for the rest lines. First take $\widehat{N}_0$ (it has all zeroes in its $0^{th}$ row and column) and fill its $0^{th}$ row (and $0^{th}$ column anti-symmetrically) as follows below to form an axillary matrix $H_0$. Make the matrix $T\widehat{N}_0 T$ where $T$ is the diagonal matrix with $diag(T) = (1, 1, (\widehat{N}_0)_{1,2}, (\widehat{N}_0)_{1,3}, \ldots, (\widehat{N}_0)_{1,n-1})$. The matrix $T\widehat{N}_0 T$ now has the 1st row filled with all +1s, (except one 0 entry on the diagonal and a zero in the 0th column) and the 1st column with all -1s, respectively. Then the 0th row should be filled with all +1s (0th column with -1s and 0 on the diagonal). Now again multiply this matrix with $T$ from both sides as before. The matrix $H_0$ that we obtained differs from the original matrix $\widehat{N}_0$ only by the 0th row (0th column) now filled with +1 and -1. The described procedure keeps the triangular structure of $H_0$, which allows the simulation of a real projection of straight oriented lines.

It is clear that in this way we can produce $2(n-1)$ matrices of type of $H_0$ from $\widehat{N}_0$ by changing the first row/column to the second one, etc., and by filling rows with -1s instead of 1s. Analogously, we can proceed with $\widehat{N}_1$ etc., to obtain in total $2n(n-1)$ matrices. Some of them may though coincide. Yet for our purpose it does not matter because we just need any of them, say, $H_0$ to form a pseudo projection matrix as:

$$[prM(\widehat{N})_i]_{j,k} = (H_0)_{i,j}(H_0)_{i,k}(\widehat{N}_i)_{j,k} \quad . (27)$$

One can make sure that $\widehat{N} = D3(prM(\widehat{N}))$ as it should be according to Eq. (16). This pseudo projection matrix also satisfies Eq. (13)

$$R(UU, prM(\widehat{N})) \equiv \widehat{0}.$$

With the help of Eq.(27) one can calculate $J_D$ and $J_M$ polynomials within this completely discrete approach. This discreteness allows one to scan all possible configurations. Moreover, it can also give the invariants for configurations, impossible for straight lines but quite possible for circles as we said before.

We have to repeat the same warning about cluster identification for the exceptional chirality matrix with $|P| = -125$: there are two clusters as shown in Table 1 and $prM(\widehat{N})$ gives $J_D$ which is the invariant of the other cluster, different from the cluster to which $Inv(P, \widehat{N})$ of the initial 3D configuration belongs. Still the clusters remain quite separable and we can distinguish them from each other.

Next we applied our approach for 7-crosses and 8-crosses to find how many topologically different configurations (rigid isotopy) can exist. As far as we know, the latter case has never been solved before, while for 7-crosses [6] reported 74 configurations. We confirm this result filling Table 2 with all 37 invariants with a positive determinant of the chiral matrix. The mirror configurations just give the negative sign to the determinant which adds additional 37 invariants to make 74 in total.

|   | $|P|$ | $J_D(|P|, i, 0.8)$ | $InvP(P)$ |
|---|---|---|---|
| 1 | 250 | -237.79522 | 8.75738 |
| 2 | 250 | -236.67421 | 8.75738 |
| 3 | 250 | -233.93737 | 8.75738 |
| 4 | 162 | -265.88901 | -1.77744 |
| 5 | 162 | -261.94412 | -1.77744 |
| 6 | 162 | -259.20728 | -1.77744 |
| 7 | 150 | -247.23958 | 6.00631 |
| 8 | 150 | -245.488 | 6.00631 |
| 9 | 102 | -262.67019 | -0.45955 |
| 10 | 102 | -258.39388 | -0.45955 |
| 11 | 90 | -251.29501 | -5.371 |
| 12 | 78 | -240.1629 | -4.03239 |
| 13 | 70 | -272.40407 | -7.9003 |
| 13 | 66 | -292.41549 | -8.02153 |
| 15 | 66 | -282.99692 | -8.02153 |
| 16 | 54 | -335.14428 | -9.34009 |
| 17 | 50 | -242.84318 | -10.07665 |
| 18 | 46 | -254.17308 | -11.50639 |
| 19 | 42 | -347.50543 | 23.2139 |
| 20 | 42 | -351.36328 | 23.2139 |
| 21 | 42 | -244.93397 | 5.91173 |
| 22 | 34 | -293.82026 | -11.78187 |
| 23 | 30 | -221.54775 | -13.32932 |
| 24 | 30 | -439.11084 | 24.24625 |
| 25 | 26 | -475.81905 | 22.03808 |
| 26 | 22 | -226.95692 | -6.36051 |
| 27 | 18 | -253.66653 | -13.90347 |
| <span style="color:red">28</span> | <span style="color:red">18</span> | <span style="color:red">-265.56832</span> | <span style="color:red">8.37689</span> |
| 29 | 18 | -236.6135 | -6.54805 |
| 30 | 18 | -233.87666 | -6.54805 |
| 31 | 14 | -592.28519 | 20.64775 |
| 32 | 10 | -300.41526 | 9.65976 |
| 33 | 10 | -222.20884 | -11.20564 |
| 34 | 6 | -169.43634 | -54.72727 |
| 35 | 2 | -271.92935 | 12.22899 |
| 36 | 2 | -305.51973 | 35.05505 |
| 37 | 2 | -304.85789 | 35.05505 |

Table 2.

We marked red in Table 2 the line 28 which shows the rigid isotopy invariant $J_D(18,1,0.8) = -265.56832$ for the mirror image (its configuration invariant $Inv(P, \widehat{N}) = -1.43470997528$ with $|P| = 18$) of the exceptional configuration of 7*-knot (its configuration invariant $Inv(P, \widehat{N}) = -1.215562687356$ with $|P| = -18$) presented in Fig. 1. Recall that only this rigid isotopy can allow 7 equal round cylinders to be in mutually touching [2].

For 8-crosses we give the complete list of rigid isotopy invariants $J_D$ in Appendix 4. One can see that we continued the series of topologically different configurations (rigid isotopy): 6-cross: 19 configurations; 7-cross: 74 configurations; 8-cross: 506 configurations. The latter result is novel.

**Conclusions**

Manifestly 3D approach to the rigid isotopy of $n$ lines in 3D reveals several important points:

- Quantization of configurations of lines allows an elementary and completely 3D description of configurations of straight lines. The number of all configurations is finite. For example, the total number of all possible configurations of 6 lines in 3D is 11618 as one can get from Table 1.
- A connection rule between adjacent configurations combines them into groupoids and distinguishes the rigid isotopy of configurations by their belonging to different groupoids.
- Quantization of configurations allows establishing connection between 3D configurations and their 2D projection diagrams.
- The tools of knot theory applied to 2D projections of line configurations lead to the same topological results as our 3D approach. A novel polynomial introduced helps to distinguish details of entanglement of lines.
- The configuration of lines --- the $n$-cross --- is naturally considered as a whole entity which is inherently fermion-like and is always a complete $n$-link of $n$ unknots in a topological sense.
- We confirmed known results for 6 and 7 lines and found that the number of topologically different rigid isotopy configurations for 8 lines is 506.

Currently, our 3D quantization of geometry and topology of lines is at the baby stage and much work is ahead.

# Appendix 1

Obtaining the Ring matrix from the direction matrix (Eq. (6)).

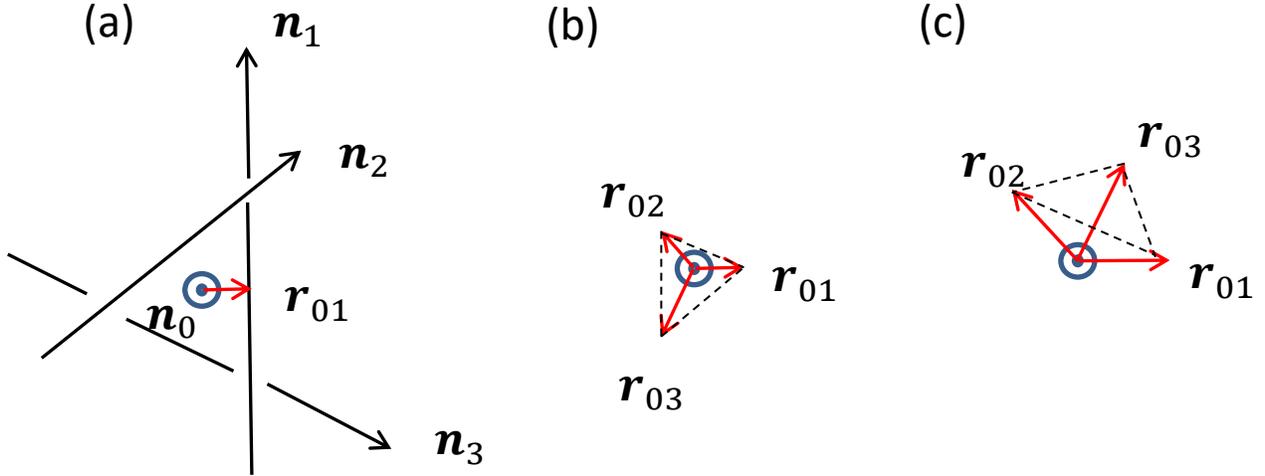

Fig. A1.1.

Consider four lines: in Fig. A1.1a the unit vector $n_0$ of line 0, directed towards the reader, is in the encircled origin, the unit vectors of three other lines are $n_1$, $n_2$, and $n_3$. Vectors directed from the origin perpendicular to three lines are $r_{01}$, $r_{02}$, and $r_{03}$ and can be defined as

$r_{0i} = P_{0i}[n_0 \times n_i]$ .  (A1.1)

In Fig. A1.1a line 0 is encaged by the three other lines. Now we can introduce an index $\mathcal{I}$ that characterizes the encaging:

$$\mathcal{I}_{0;1,2,3} = \frac{|[r_{01} \times r_{02}] + [r_{02} \times r_{03}] + [r_{03} \times r_{01}]|}{|[r_{01} \times r_{02}]| + |[r_{02} \times r_{03}]| + |[r_{03} \times r_{01}]|} .$$  (A1.2)

As it is seen from Fig. A1.1b, this index is 1 when the line 0 is encaged, because the total area covered with three triangles coincide with the area with the dashed triangle, while for the case in Fig. A1c the index is always less than one. Using the identity

$[n_0 \times n_i] \times [n_0 \times n_j] = n_0(n_i, n_j, n_0)$  (A1.3)

and the definition (A1.1) we further rewrite Eq. (A1.2) as

$$\mathcal{I}_{0;1,2,3} = \frac{|P_{10}P_{20}(n_1,n_2,n_0) + P_{20}P_{30}(n_2,n_3,n_0) + P_{30}P_{10}(n_3,n_1,n_0)|}{|(n_1,n_2,n_0)| + |(n_2,n_3,n_0)| + |(n_3,n_1,n_0)|} .$$  (A1.4)

Now Eq. (A1.4) gives the value 1 when line zero is encaged and a value less than 1 if not. Still one would prefer to have an indicator of encaging that would have the other value exactly zero. We found that it is possible to modify Eq. (A1.4) exactly in such a way by noticing that if the terms in the nominator are of the same sign, the index will be definitely 1, otherwise 0:

$$\mathcal{I}_{0;1,2,3} =$$

$$\frac{1}{8}\{[P_{10}P_{20}sign(\boldsymbol{n}_1,\boldsymbol{n}_2,\boldsymbol{n}_0)+1][P_{20}P_{30}sign(\boldsymbol{n}_2,\boldsymbol{n}_3,\boldsymbol{n}_0)+1][P_{30}P_{10}sign(\boldsymbol{n}_3,\boldsymbol{n}_1,\boldsymbol{n}_0)+1] -$$
$$[P_{10}P_{20}sign(\boldsymbol{n}_1,\boldsymbol{n}_2,\boldsymbol{n}_0)-1][P_{20}P_{30}sign(\boldsymbol{n}_2,\boldsymbol{n}_3,\boldsymbol{n}_0)-1][P_{30}P_{10}sign(\boldsymbol{n}_3,\boldsymbol{n}_1,\boldsymbol{n}_0)-1]\} .$$
(A1.5)

One may recognize the elements of the direction matrix appeared in Eq. (A1.5) $(\widehat{N}_i)_{j,k}$ so that Eq. (A1.5) takes a more general form:

$$\mathcal{I}_{l;i,j,k} = \frac{1}{8}\{[P_{il}P_{jl}(\widehat{N}_l)_{i,j}+1][P_{jl}P_{kl}(\widehat{N}_l)_{j,k}+1][P_{kl}P_{il}(\widehat{N}_l)_{k,i}+1] - [P_{il}P_{jl}(\widehat{N}_l)_{i,j}-1][P_{jl}P_{kl}(\widehat{N}_l)_{j,k}-1][P_{kl}P_{il}(\widehat{N}_l)_{k,i}-1]\} . \quad (A1.6)$$

Now it is clear that the Ring matrix entry can be obtained through $\mathcal{I}_{l;i,j,k}$ as

$$R_{li} = \frac{1}{2}\Sigma_{j,k}\mathcal{I}_{l;i,j,k}(1-\delta_{l,i}), \quad (A1.7)$$

because according to its definition, the Ring matrix entry $R_{li}$ is the number of times when $i$th line participates in encircling the $l$th line in different triangles made with the help of $j$th and $k$th lines. Factor one-half is introduced not to count twice because of the symmetry and the diagonal is made zero with the Kronecker delta.

Using Eq. (A1.6) and (A1.7) after some algebra we obtain Eq. (6) of the main text.

**Appendix 2**

| | | |
|---|---|---|
| 0 | -2.37579 | 8, 4, 6, 21 |
| 1 | -2.37797 | 4 |
| 2 | -2.37964 | 11 |
| 3 | -2.35219 | 8,18,4 |
| 4 | -2.22147 | 1,0,18,6,3 |
| 5 | -2.21013 | 10,6,8,32,40 |
| 6 | -2.13312 | 0,4,5 |
| 7 | -2.09907 | 12,41,13,36 |
| 8 | -2.09829 | 0,31,3,19,5 |
| 9 | -2.06866 | 12,10,40,32 |
| 10 | -2.055 | 5,9,39,34 |
| 11 | -2.01438 | 2,22,16 |
| 12 | -1.96632 | 48,38,9,39,7 |
| 13 | -1.89822 | 38,7 |
| 14 | -1.8718 | 48,24 |
| 15 | -1.85554 | 35,20,17,25 |
| 16 | -1.84469 | 11,42 |
| 17 | -1.83636 | 30,15,35 |
| 18 | -1.80519 | 3,45,4 |
| 19 | -1.75345 | 46,32,8,21 |
| 20 | -1.7421 | 15,43,28 |
| 21 | <span style="color:red">-1.74039</span> | 0,31,34,19 |
| 22 | -1.73732 | 11,33,43 |
| 23 | -1.69622 | 33,26,28,43,47 |
| 24 | -1.65429 | 39,14,48,29 |
| 25 | -1.64356 | 15,37,42 |
| 26 | -1.64029 | 31,23,27 |
| 27 | -1.58079 | 31,28,44,45 |
| 28 | -1.58054 | 20,27,23,35 |
| 29 | -1.53648 | 43,47,33,31 |
| 30 | -1.52182 | 37,17 |
| 31 | -1.45052 | 8,47,26,21,29 |
| 32 | -1.44442 | 5,19,9,38,34 |
| 33 | -1.43336 | 22,23,29 |
| 34 | -1.40962 | 21,32,10 |
| 35 | -1.4084 | 17,15,28,44 |
| 36 | -1.40598 | 7,48,39 |
| 37 | -1.37042 | 30,25 |
| 38 | -1.35073 | 12,13,46,32 |
| 39 | -1.27873 | 47,36,24,10,12 |
| 40 | 0.1123 | 9,5 |
| 41 | 0.18389 | 7 |
| 42 | 0.84907 | 25,16 |
| 43 | 0.91978 | 29,20,23,22 |
| 44 | 1.34007 | 35,27 |
| 45 | 2.28626 | 27,18 |
| 46 | 2.32991 | 19,38 |
| 47 | 3.39671 | 48,39,23,31,29 |
| 48 | -1.73412 | 47,14,24,36,12 |

Table A2.1. All 49 invariants of the cluster 27** from Table 1. The invariant at number 21 marked red is obtained for the corresponding line configuration of [3]. In the third column the connected close neighbors are given to make sure that the groupoid in fully connected.

## Appendix 3

Let us show the calculations of $J_M(2, 0, a)$ for a 3-cross.

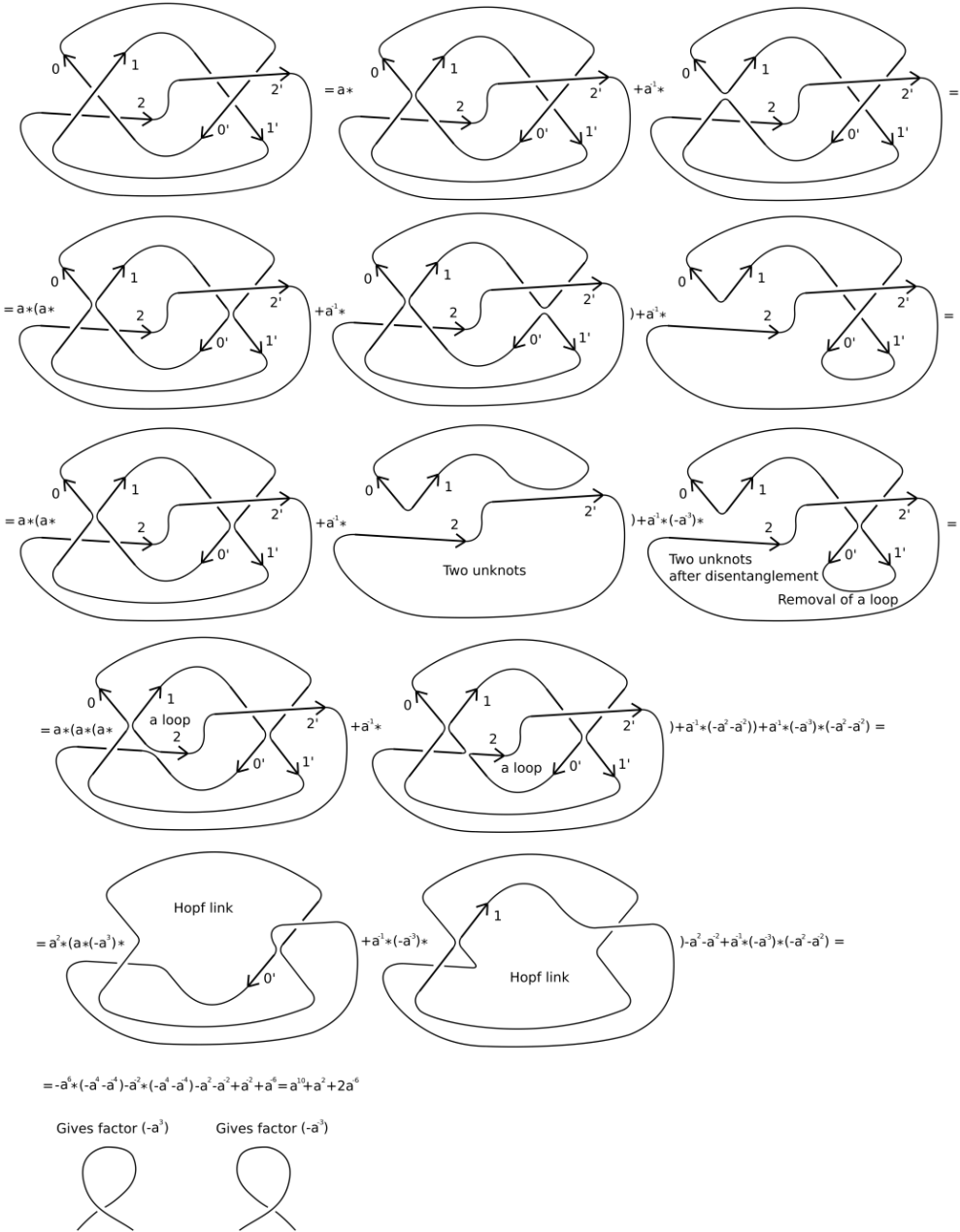

The result is $J_M(2,0,a) = a^{10} + a^2 + 2a^{-6}$.

Next let us calculate $J_D(-2,0,a^{-1})$ for a 3-cross (it is the mirror image of the above 3-cross).

$= -a^{-1} - a^3 - a^{-1} + a^3 - a^{-1} - a^{-5} = -2a^3 - a^{-1} - a^{-5}$

The result is $J_D(-2,0,a^{-1}) = -2a^3 - a^{-1} - a^{-5}$. Notice that $J_D(2,0,a) = -2a^{-3} - a^1 - a^5$ and $-J_D(2,0,a^2) = J_M(2,0,a)$.

# Appendix 4

| | $\lvert P \rvert$ | $J_D(\lvert P \rvert, i, 0.8)$ | $InvP(P)$ |
|---|---|---|---|
| 1 | -495 | -602.7108580028 | 38.47877 |
| 2 | -495 | -613.5641767248 | -6.93979 |
| 3 | -495 | -605.0509190290 | -6.93979 |
| 4 | -495 | -609.7339582672 | -6.93979 |
| 5 | -495 | -593.6177177638 | -6.93979 |
| 6 | -495 | -917.7737832429 | -6.93979 |
| 7 | -495 | -929.2069845080 | 38.47877 |
| 8 | -495 | -921.7970450096 | 38.47877 |
| 9 | -495 | -889.8606942165 | 38.47877 |
| 10 | -495 | -904.4718180932 | 38.47877 |
| 11 | -375 | -589.7584848208 | 36.96626 |
| 12 | -375 | -578.3252835556 | -6.4565 |
| 13 | -375 | -584.9599784679 | -6.4565 |
| 14 | -375 | -589.6430177061 | -6.4565 |
| 15 | -375 | -1099.6391071643 | -6.4565 |
| 16 | -375 | -1088.2059058992 | 36.96626 |
| 17 | -375 | -1096.7038528152 | 36.96626 |
| 18 | -375 | -1068.7907637889 | 36.96626 |
| 19 | -351 | -774.3566791687 | -14.31083 |
| 20 | -351 | -779.5476899588 | -14.31083 |
| 21 | -351 | -781.1068411956 | -14.31083 |
| 22 | -351 | -777.0191486029 | -14.31083 |
| 23 | -351 | -711.2006846632 | 15.77123 |
| 24 | -351 | -716.3916954333 | 15.77123 |
| 25 | -351 | -732.8715831944 | 15.77123 |
| 26 | -351 | -733.8537690739 | 15.77123 |
| 27 | -295 | -659.0158692321 | -37.71272 |
| 28 | -295 | -656.6874982457 | -37.71272 |
| 29 | -295 | -671.1854722813 | -37.71272 |
| 30 | -295 | -649.3702494359 | -37.71272 |
| 31 | -295 | -662.6722145855 | -37.71272 |
| 32 | -295 | -851.9559953159 | 26.24782 |
| 33 | -295 | -846.0176251509 | 26.24782 |
| 34 | -295 | -834.0916183391 | 26.24782 |
| 35 | -295 | -843.4427376202 | 26.24782 |
| 36 | -295 | -830.1407724706 | 26.24782 |
| 37 | -279 | -681.6609138345 | -18.16848 |
| 38 | -279 | -698.1408015956 | -18.16848 |
| 39 | -279 | -670.2277125693 | -18.16848 |
| 40 | -279 | -691.0120331155 | -18.16848 |
| 41 | -279 | -683.5296777189 | 16.93263 |
| 42 | -279 | -876.4999014962 | 16.93263 |
| 43 | -279 | -873.5799579268 | 16.93263 |
| 44 | -279 | -885.0131591920 | 16.93263 |
| 45 | -279 | -880.3301199538 | 16.93263 |
| 46 | -279 | -879.5646746267 | -18.16848 |
| 47 | -255 | -609.1678009765 | -36.78619 |
| 48 | -255 | -614.3588117665 | -36.78619 |
| 49 | -255 | -616.4850497862 | 31.4695 |
| 50 | -255 | -967.6711175166 | -36.78619 |
| 51 | -255 | -954.9977513298 | 31.4695 |
| 52 | -255 | -972.8621283067 | 31.4695 |
| 53 | -231 | -704.8412334131 | -13.47155 |
| 54 | -231 | -694.4745226126 | 27.17499 |
| 55 | -231 | -910.3829480063 | 27.17499 |
| 56 | -231 | -912.2615107898 | -13.47155 |
| 57 | -199 | -729.2602077817 | -18.05874 |
| 58 | -199 | -721.1492534222 | 21.6594 |
| 59 | -199 | -703.2848764454 | 21.6594 |
| 60 | -199 | -905.8706223253 | 21.6594 |
| 61 | -199 | -898.5533735156 | -18.05874 |
| 62 | -199 | -895.2311266100 | -18.05874 |
| 63 | -175 | -630.4647961184 | -16.67047 |
| 64 | -175 | -627.1425492127 | -16.67047 |
| 65 | -175 | -663.5205076632 | -37.71254 |
| 66 | -175 | -645.6561306864 | -37.71254 |
| 67 | -175 | -658.9580958360 | 28.62347 |
| 68 | -175 | -995.9155375737 | 28.62347 |
| 69 | -175 | -987.4022798779 | 25.8705 |
| 70 | -175 | -988.5982887640 | 25.8705 |
| 71 | -175 | -998.3814140160 | 25.8705 |
| 72 | -175 | -1006.9923683755 | -16.67047 |
| 73 | -175 | -666.1296645716 | 19.47846 |
| 74 | -175 | -1066.8918143526 | -22.84048 |
| 75 | -159 | -680.9468049627 | 23.82026 |
| 76 | -159 | -663.0824279859 | 17.09527 |
| 77 | -159 | -843.6865721531 | -21.15604 |
| 78 | -159 | -810.5824723367 | 17.09527 |
| 79 | -159 | -890.3626789175 | -16.84516 |
| 80 | -159 | -880.4886288892 | 23.82026 |
| 81 | -159 | -1026.7767553161 | -16.84516 |
| 82 | -159 | -1019.4595065064 | -21.15604 |
| 83 | -151 | -614.7389837133 | -21.37839 |
| 84 | -151 | -1172.1663802374 | 22.87575 |
| 85 | -135 | -649.0711276886 | 18.07059 |
| 86 | -135 | -659.0508459326 | -25.22167 |
| 87 | -135 | -645.7488807830 | -23.24844 |
| 88 | -135 | -654.2621384787 | -23.24844 |
| 89 | -135 | -658.4222469697 | -20.35343 |
| 90 | -135 | -1119.1627363951 | 19.62765 |
| 91 | -135 | -852.4233730461 | 19.62765 |
| 92 | -135 | -847.2323622560 | -20.35343 |
| 93 | -135 | -860.5343274057 | 19.62765 |
| 94 | -135 | -858.4080893860 | 26.27575 |
| 95 | -135 | -852.0210697099 | -23.24844 |
| 96 | -135 | -704.9463224030 | 26.27575 |
| 97 | -135 | -610.6446930170 | -23.24844 |
| 98 | -135 | -605.4536822269 | 19.62765 |
| 99 | -135 | -1111.6306597743 | 21.46702 |
| 100 | -135 | -1116.8216705644 | -21.67624 |
| 101 | -135 | -629.2472672858 | 19.62765 |
| 102 | -135 | -1207.5869944281 | -23.24844 |
| 103 | -119 | -648.6247806801 | -22.67008 |
| 104 | -119 | -1280.6281828493 | 15.60047 |
| 105 | -111 | -730.3365201657 | 14.87095 |
| 106 | -111 | -690.9902298742 | 23.20487 |
| 107 | -111 | -715.0967973261 | 14.87095 |
| 108 | -111 | -1067.7920953956 | -22.71231 |
| 109 | -111 | -1065.2352941771 | -24.55959 |
| 110 | -111 | -1081.3515346804 | 14.87095 |
| 111 | -111 | -581.8489479640 | -22.71231 |
| 112 | -111 | -1336.6578894325 | -22.71231 |
| 113 | -111 | -704.1840095114 | -19.91999 |
| 114 | -111 | -1138.8058667363 | 17.26775 |
| 115 | -103 | -755.1012166903 | -24.528 |
| 116 | -103 | -794.4475069818 | 15.72441 |
| 117 | -103 | -944.2912934566 | 15.72441 |
| 118 | -103 | -960.4075339600 | -24.528 |
| 119 | -87 | -928.1821865364 | 16.13766 |
| 120 | -87 | -843.8197361988 | -27.2955 |
| 121 | -79 | -611.7423833943 | 28.54075 |
| 122 | -79 | -1421.3233629248 | 15.74075 |
| 123 | -79 | -639.5562803150 | 4.16907 |
| 124 | -79 | -647.6672346746 | -26.38317 |
| 125 | -79 | -1002.0153879451 | 4.16907 |
| 126 | -79 | -998.9314140394 | 28.54075 |
| 127 | -71 | -576.1844769423 | 13.07421 |
| 128 | -71 | -1642.9207918514 | -25.60311 |
| 129 | -63 | -622.9505294896 | 32.58225 |
| 130 | -63 | -1448.6870008527 | -27.70136 |
| 131 | -63 | -699.9868744112 | 14.70753 |
| 132 | -63 | -715.5801841673 | 13.03651 |
| 133 | -63 | -708.0973287707 | -18.269 |
| 134 | -63 | -628.4034500034 | 17.18544 |
| 135 | -63 | -634.5766467130 | -18.269 |
| 136 | -63 | -633.5944608335 | 5.07086 |
| 137 | -63 | -948.6037244325 | 53.30521 |
| 138 | -63 | -953.7947352226 | 5.07086 |
| 139 | -63 | -944.5160318398 | 5.07086 |
| 140 | -63 | -786.7757709446 | 53.30521 |
| 141 | -63 | -783.4535240389 | 17.18544 |
| 142 | -63 | -786.5182968094 | 17.18544 |
| 143 | -63 | -1253.6129464023 | 13.03651 |
| 144 | -63 | -1248.4219356122 | 32.58225 |
| 145 | -63 | -578.1023011464 | 32.58225 |
| 146 | -63 | -583.2933119365 | 13.03651 |
| 147 | -55 | -657.8139814670 | 51.59818 |
| 148 | -55 | -661.1362283726 | 15.76878 |
| 149 | -55 | -881.8791852630 | -18.22209 |
| 150 | -55 | -889.9901396225 | 51.59818 |
| 151 | -55 | -1045.9288422676 | 30.09419 |
| 152 | -55 | -1035.7676565521 | 51.59818 |
| 153 | -55 | -1057.3620435327 | -29.59457 |
| 154 | -55 | -1040.7378314775 | 51.59818 |
| 155 | -55 | -627.7911502873 | -18.22209 |
| 156 | -55 | -622.6001394972 | -27.8302 |
| 157 | -55 | -611.1669382221 | -27.8302 |
| 158 | -55 | -641.6151442616 | -18.22209 |
| 159 | -55 | -786.8535965071 | -18.22209 |
| 160 | -55 | -1046.9538605081 | 30.09419 |
| 161 | -47 | -736.2050213577 | -29.12583 |
| 162 | -47 | -1191.7392092008 | 12.97563 |
| 163 | -39 | -554.4975668836 | -22.52487 |
| 164 | -39 | -1473.4491652704 | -28.84552 |
| 165 | -39 | -602.6738830293 | 32.11652 |
| 166 | -39 | -1142.7105680206 | 10.37856 |
| 167 | -39 | -520.5373810637 | 5.26753 |
| 168 | -39 | -2153.6581299942 | 46.78531 |
| 169 | -31 | -765.8934595175 | 21.77214 |
| 170 | -31 | -745.1091389712 | 44.71527 |
| 171 | -31 | -645.4104119676 | -22.4469 |
| 172 | -31 | -1055.1641514844 | -0.963 |
| 173 | -31 | -837.0890462947 | -22.4469 |
| 174 | -31 | -831.6405613695 | 44.71527 |
| 175 | -31 | -1234.8446750656 | 22.37549 |
| 176 | -31 | -1246.2778763308 | 1.69773 |
| 177 | -31 | -596.8705724832 | 1.69773 |
| 178 | -31 | -585.4373712180 | 22.37549 |
| 179 | -23 | -607.9368567591 | -30.89264 |
| 180 | -23 | -1729.4670583995 | 8.39344 |
| 181 | -15 | -798.9111127582 | 19.96149 |
| 182 | -15 | -937.9301213923 | 37.66027 |
| 183 | -15 | -921.7166520574 | -26.54396 |
| 184 | -15 | -742.9618324982 | -2.71931 |
| 185 | -15 | -1730.4671483877 | 19.59178 |
| 186 | -15 | -533.4479440040 | 0.39626 |
| 187 | -7 | -682.2928132533 | 46.54545 |
| 188 | -7 | -690.4037676128 | -78.76923 |
| 189 | -7 | -749.0477608044 | 23.49134 |
| 190 | -7 | -747.0863062313 | 23.49134 |
| 191 | -7 | -748.4430849348 | 26.54146 |
| 192 | -7 | -750.8685252416 | 23.49134 |
| 193 | -7 | -723.2704106157 | 26.54146 |
| 194 | -7 | -751.8133312878 | 2.66161 |
| 195 | -7 | -911.1467245031 | 34.32786 |
| 196 | -7 | -914.4689714088 | 34.32786 |
| 197 | -7 | -392.3615445342 | 2.66161 |
| 198 | -7 | -4347.0960429704 | 34.32786 |
| 199 | 1 | -901.9654727023 | 23.00485 |
| 200 | 1 | -1120.0405778919 | 30.03321 |
| 201 | 1 | -730.2169261505 | 2.49266 |
| 202 | 1 | -689.9828251712 | 30.03321 |
| 203 | 1 | -836.8226786406 | -0.44929 |
| 204 | 1 | -839.5875450068 | 19.804 |
| 205 | 1 | -671.9838002241 | 24.24694 |
| 206 | 1 | -771.6825232278 | 24.24694 |
| 207 | 9 | -949.8895855995 | 0.71026 |
| 208 | 9 | -808.9161127526 | 2.11661 |
| 209 | 9 | -1002.0145908984 | 24.64729 |
| 210 | 9 | -951.2578451001 | 2.17457 |
| 211 | 9 | -955.5040498440 | 24.53407 |
| 212 | 9 | -766.1115241304 | 28.17424 |
| 213 | 9 | -783.5022356633 | 56.59498 |
| 214 | 9 | -731.5353905919 | 26.91336 |
| 215 | 9 | -699.9901207636 | 56.59498 |
| 216 | 9 | -695.4037858430 | 56.59498 |
| 217 | 9 | -847.2191549036 | 21.99287 |
| 218 | 9 | -1063.6444986125 | 2.11661 |
| 219 | 9 | -1065.6059531856 | 23.62488 |
| 220 | 9 | -1070.7969639757 | 56.59498 |
| 221 | 9 | -1066.4768602784 | 24.64729 |
| 222 | 9 | -696.8010968258 | 2.11661 |
| 223 | 9 | -602.9230285792 | 2.11661 |
| 224 | 9 | -597.7320177891 | -2.60287 |
| 225 | 9 | -591.2318482720 | 4.55835 |
| 226 | 9 | -571.9821931624 | 4.55835 |
| 227 | 17 | -853.7448199798 | 25.82892 |
| 228 | 17 | -846.4275711701 | 2.602 |
| 229 | 17 | -814.0356370277 | 2.602 |
| 230 | 17 | -812.7803061009 | 32.24894 |
| 231 | 17 | -690.3533400312 | -53.65149 |
| 232 | 17 | -708.6208092613 | 28.42325 |
| 233 | 17 | -778.3940699330 | 24.61858 |
| 234 | 17 | -760.5296929561 | 25.82892 |
| 235 | 17 | -450.3986451200 | 24.61858 |
| 236 | 17 | -3045.7209061680 | 32.24894 |
| 237 | 25 | -1007.2400385871 | 26.54738 |
| 238 | 25 | -1261.7764630305 | 2.86272 |
| 239 | 25 | -979.3453978163 | -2.18605 |
| 240 | 25 | -813.4199064406 | 18.89379 |
| 241 | 25 | -708.8661975098 | 20.35611 |
| 242 | 25 | -720.2918423695 | 23.15325 |
| 243 | 25 | -629.7730608613 | 0.41694 |
| 244 | 25 | -967.1098439352 | -11.13043 |
| 245 | 33 | -959.8945855940 | -50.13708 |
| 246 | 33 | -957.9931310209 | 0.93513 |
| 247 | 33 | -881.1084579077 | 36.37988 |
| 248 | 33 | -610.6124193648 | 15.30809 |
| 249 | 33 | -622.3035996720 | 15.30809 |
| 250 | 33 | -771.7378390514 | 28.63547 |
| 251 | 33 | -493.1859999832 | 31.11657 |
| 252 | 33 | -2290.4494058321 | 36.37988 |
| 253 | 41 | -515.9430601125 | 11.14174 |
| 254 | 41 | -1944.2035168844 | 20.13391 |
| 255 | 41 | -885.0131591920 | 0.05023 |
| 256 | 41 | -880.4355813596 | 3.33234 |
| 257 | 41 | -881.6909128863 | 11.14174 |
| 258 | 41 | -767.3951998518 | 20.13391 |
| 259 | 41 | -760.5859225939 | 26.09961 |
| 260 | 41 | -960.4936779976 | 26.09961 |
| 261 | 41 | -951.8953289141 | 26.09961 |
| 262 | 41 | -944.5780801044 | 11.14174 |
| 263 | 41 | -776.7156165285 | 3.33234 |
| 264 | 41 | -728.8709178600 | 24.69132 |
| 265 | 41 | -711.0065408831 | 24.69132 |
| 266 | 41 | -690.0144776911 | 0.77991 |
| 267 | 41 | -846.2846972158 | 0.77991 |
| 268 | 41 | -805.6982420028 | 20.94075 |
| 269 | 41 | -678.3386669288 | -50.5016 |
| 270 | 41 | -660.0711976987 | 33.89381 |
| 271 | 41 | -670.2277125693 | 25.09239 |
| 272 | 49 | -1255.8660184316 | 1.26663 |
| 273 | 49 | -717.0777432861 | -9.58863 |
| 274 | 49 | -1097.2849538355 | -6.39887 |
| 275 | 49 | -662.6630966106 | 25.309 |
| 276 | 57 | -600.4592511129 | -46.63787 |
| 277 | 57 | -575.2776253264 | -46.63787 |
| 278 | 57 | -592.9615849366 | 32.82887 |
| 279 | 57 | -1450.2496414593 | 32.82887 |
| 280 | 57 | -1430.5067016990 | 32.82887 |
| 281 | 57 | -1425.0680156727 | -46.63787 |
| 282 | 65 | -782.7232511909 | 25.44977 |
| 283 | 65 | -783.5266629841 | 25.44977 |
| 284 | 65 | -785.6529010037 | 25.44977 |
| 285 | 65 | -775.1781887348 | 25.44977 |
| 286 | 65 | -777.5419466442 | -18.60537 |
| 287 | 65 | -777.1396433079 | -18.60537 |
| 288 | 65 | -839.5875450068 | -46.82774 |
| 289 | 65 | -836.2652983012 | -18.60537 |
| 290 | 65 | -730.2169261505 | -46.82774 |
| 291 | 65 | -738.3278805101 | 25.44977 |
| 292 | 65 | -706.9749588298 | -18.60537 |
| 293 | 65 | -695.2837785225 | 25.44977 |
| 294 | 65 | -706.3463598669 | -18.60537 |
| 295 | 65 | -693.6729936802 | 3.8456 |
| 296 | 65 | -696.9952405859 | 24.24406 |
| 297 | 65 | -677.8034391949 | -18.60537 |
| 298 | 65 | -649.8326497115 | 24.24406 |
| 299 | 65 | -607.6927134220 | 33.72698 |
| 300 | 65 | -632.8743392085 | 33.72698 |
| 301 | 65 | -1224.0225341268 | 33.72698 |
| 302 | 65 | -1187.7527631137 | 3.8456 |
| 303 | 65 | -1198.8409083402 | -46.82774 |
| 304 | 73 | -905.8688668301 | -5.51396 |
| 305 | 81 | -808.8446263376 | 26.09279 |
| 306 | 81 | -811.9093990081 | 26.12711 |
| 307 | 81 | -814.0356370277 | 26.09279 |
| 308 | 81 | -870.4487627492 | 6.74517 |
| 309 | 81 | -747.6432234500 | 6.74517 |
| 310 | 81 | -713.8118200514 | 6.74517 |
| 311 | 81 | -708.6208092613 | -4.79144 |
| 312 | 81 | -721.2941754480 | -7.52447 |
| 313 | 81 | -819.0575717811 | 26.09279 |
| 314 | 81 | -996.6351491424 | -7.01963 |
| 315 | 81 | -1540.81071590190 | -3.47278 |
| 316 | 81 | -604.3574104844 | 3.61219 |
| 317 | 89 | -931.0378657808 | 4.23995 |
| 318 | 89 | -795.5673363473 | 24.80223 |
| 319 | 89 | -896.4617322422 | 9.67097 |
| 320 | 89 | -677.9493770071 | 22.72589 |
| 321 | 89 | -795.5673363473 | -5.0722 |
| 322 | 105 | -827.4771635296 | 27.61329 |
| 323 | 105 | -822.2861527395 | 27.61329 |
| 324 | 105 | -814.1751983800 | -2.14387 |
| 325 | 105 | -801.5018321932 | 25.42344 |
| 326 | 105 | -815.1336873763 | 9.76427 |
| 327 | 105 | -1106.8545898520 | 25.42344 |
| 328 | 105 | -1058.6540533279 | 27.61329 |
| 329 | 105 | -785.5715094519 | 9.76427 |
| 330 | 105 | -772.8981432652 | -10.01068 |
| 331 | 105 | -772.8981432652 | -10.01068 |
| 332 | 105 | -770.7719052456 | -5.31879 |
| 333 | 105 | -715.7805287974 | -10.01068 |
| 334 | 105 | -723.0238786537 | -10.01068 |
| 335 | 105 | -614.7349027790 | 27.61329 |
| 336 | 105 | -624.8914176496 | -2.14387 |
| 337 | 105 | -630.0824284397 | -10.01068 |
| 338 | 105 | -633.4046753454 | 27.61329 |
| 339 | 105 | -635.5309133650 | -5.31879 |
| 340 | 121 | -1072.1164649946 | -3.39343 |
| 341 | 121 | -746.8652266261 | -3.6304 |
| 342 | 129 | -1196.4260528338 | -2.73514 |
| 343 | 129 | -673.7932695209 | 8.17286 |
| 344 | 129 | -723.3749911922 | -4.76284 |
| 345 | 129 | -871.0486521995 | 27.97166 |
| 346 | 137 | -814.4190561660 | 2.10549 |
| 347 | 153 | -747.2598043117 | 9.1702 |
| 348 | 153 | -741.8113193864 | -2.75489 |
| 349 | 153 | -743.7727739595 | 9.1702 |
| 350 | 153 | -887.2921991296 | 9.1702 |
| 351 | 153 | -825.8136205491 | 25.03198 |
| 352 | 153 | -815.4992266270 | -2.75489 |
| 353 | 153 | -797.1260166980 | -2.75489 |
| 354 | 153 | -785.4348363908 | 25.03198 |
| 355 | 153 | -764.6505158445 | 25.03198 |
| 356 | 161 | -856.1564393788 | -2.10747 |
| 357 | 161 | -953.3506753772 | -2.37626 |
| 358 | 161 | -1014.8292539576 | 0.52129 |
| 359 | 161 | -759.6176061833 | -0.02155 |
| 360 | 161 | -749.3032122611 | -2.37626 |
| 361 | 161 | -725.8734899043 | 10.12996 |
| 362 | 161 | -904.4639714144 | 21.43552 |
| 363 | 161 | -756.7903104071 | -2.10747 |
| 364 | 169 | -873.7925391601 | -1.82259 |
| 365 | 177 | -737.0391310520 | 5.83341 |
| 366 | 177 | -734.9128930323 | 1.0939 |
| 367 | 177 | -889.9087480707 | 5.83341 |
| 368 | 177 | -877.2353818840 | 1.0939 |
| 369 | 185 | -745.2540609971 | -1.6995 |
| 370 | 185 | -750.4450717872 | 6.96307 |
| 371 | 185 | -837.0882492481 | 6.96307 |
| 372 | 201 | -995.1749379031 | 4.57016 |
| 373 | 201 | -690.8360301527 | 3.50145 |
| 374 | 209 | -822.1412307136 | 20.78805 |
| 375 | 209 | -820.0145926940 | 20.78805 |
| 376 | 209 | -729.3516471124 | 0.24224 |
| 377 | 209 | -734.0346863506 | 0.24224 |
| 378 | 209 | -923.0563568272 | 2.3866 |
| 379 | 209 | -895.1432678009 | 2.3866 |
| 380 | 209 | -729.3997691340 | 19.79148 |
| 381 | 209 | -716.7264029472 | 19.79148 |
| 382 | 225 | -764.8959040930 | 1.88019 |
| 383 | 225 | -912.5695651003 | 2.36139 |
| 384 | 233 | -802.4672715709 | 4.23458 |
| 385 | 233 | -777.6596892579 | 4.23458 |
| 386 | 233 | -810.2558175224 | 4.23458 |
| 387 | 233 | -782.3427284961 | 4.23458 |
| 388 | 249 | -817.7633426868 | 4.54001 |
| 389 | 249 | -845.6764317131 | 5.82393 |
| 390 | 249 | -796.8510120686 | 4.54001 |
| 391 | 249 | -792.1679728304 | 5.82393 |
| 392 | 281 | -795.5673363473 | 6.09731 |
| 393 | 297 | -903.3700100567 | -5.04804 |
| 394 | 297 | -933.7061696040 | -5.04804 |
| 395 | 297 | -902.7653343871 | 3.18651 |
| 396 | 297 | -916.0433762434 | -5.04804 |
| 397 | 297 | -907.9563449772 | -5.04804 |
| 398 | 297 | -912.1861459620 | 3.21321 |
| 399 | 297 | -659.0976571883 | 18.60595 |
| 400 | 297 | -652.1022857201 | -5.04804 |
| 401 | 297 | -653.0470917663 | 3.21321 |
| 402 | 297 | -649.9760477005 | 3.21321 |
| 403 | 297 | -647.8560809762 | 3.21321 |
| 404 | 297 | -652.1761846735 | 3.21321 |
| 405 | 305 | -825.3264981892 | 2.12544 |
| 406 | 305 | -820.1354873991 | -18.65896 |
| 407 | 305 | -688.0961237994 | 4.31375 |
| 408 | 305 | -696.2070783190 | -18.65896 |
| 409 | 305 | -699.5293250646 | -18.65896 |
| 410 | 305 | -819.0126174582 | -8.78978 |
| 411 | 305 | -822.3348643638 | -8.78978 |
| 412 | 305 | -810.9016630987 | 4.31375 |
| 413 | 305 | -697.3299480998 | 2.12544 |
| 414 | 305 | -702.5209588899 | 4.31375 |
| 415 | 329 | -796.1756018500 | 5.12966 |
| 416 | 329 | -783.5022356633 | 16.07901 |
| 417 | 329 | -766.1115241304 | 16.07901 |
| 418 | 329 | -763.9852861108 | 5.12966 |
| 419 | 345 | -733.5217379846 | 8.0474 |
| 420 | 345 | -744.9649392498 | 8.0474 |
| 421 | 345 | -739.5064543245 | 8.0474 |
| 422 | 345 | -771.6967700637 | -15.56659 |
| 423 | 345 | -750.9124495174 | -15.56659 |
| 424 | 345 | -762.3456507826 | -15.56659 |
| 425 | 369 | -701.8341901265 | -2.59044 |
| 426 | 369 | -700.4733977940 | -2.59044 |
| 427 | 369 | -698.5119432209 | 2.81364 |
| 428 | 369 | -872.5420514520 | 2.81364 |
| 429 | 369 | -846.7922268252 | -2.59044 |
| 430 | 369 | -853.8727788857 | -2.59044 |
| 431 | 369 | -828.9278498484 | -2.59044 |
| 432 | 369 | -826.8016118288 | 2.81364 |
| 433 | 369 | -703.1251654367 | 2.81364 |
| 434 | 369 | -711.8495527539 | 2.81364 |
| 435 | 369 | -718.7956767359 | 8.32472 |
| 436 | 369 | -698.8050617395 | 12.99198 |
| 437 | 369 | -706.1223105492 | 8.32472 |
| 438 | 369 | -824.6397294258 | 8.32472 |
| 439 | 369 | -832.7506837853 | 12.99198 |
| 440 | 369 | -844.4418640925 | 12.99198 |
| 441 | 377 | -754.0464234464 | 13.26706 |
| 442 | 377 | -766.7197896332 | 13.26706 |
| 443 | 377 | -751.9201854268 | 13.26706 |
| 444 | 425 | -774.6918818305 | -14.76235 |
| 445 | 425 | -771.7566274814 | -14.76235 |
| 446 | 425 | -772.5024125949 | -12.1589 |
| 447 | 425 | -779.1371075072 | 7.62738 |
| 448 | 425 | -769.2967916060 | -12.1589 |
| 449 | 425 | -764.1057808159 | 7.62738 |
| 450 | 441 | -861.2024503261 | 50.32306 |
| 451 | 441 | -860.8154577696 | 12.6203 |
| 452 | 441 | -655.2737431764 | 12.6203 |
| 453 | 441 | -653.7974837292 | 50.32306 |
| 454 | 505 | -746.4414167917 | 12.9024 |
| 455 | 505 | -748.9699581476 | 12.9024 |
| 456 | 505 | -741.6527093379 | 12.9024 |
| 457 | 505 | -723.7883323610 | 12.9024 |
| 458 | 545 | -783.6417970155 | 10.38462 |
| 459 | 545 | -791.7527513751 | 10.38462 |
| 460 | 545 | -793.6215153595 | 1.25248 |
| 461 | 545 | -796.9437621651 | 10.38462 |
| 462 | 545 | -682.3821325188 | 10.38462 |
| 463 | 545 | -690.8953902145 | 1.25248 |
| 464 | 545 | -690.4930868783 | 1.25248 |
| 465 | 545 | -687.5731433088 | 1.25248 |
| 466 | 585 | -746.2044758451 | -5.09795 |
| 467 | 585 | -726.4023411783 | -5.09795 |
| 468 | 585 | -743.8644148189 | -5.09795 |
| 469 | 585 | -730.2325996359 | -5.09795 |
| 470 | 585 | -734.2558214027 | -5.09795 |
| 471 | 585 | -723.0800942727 | -5.09795 |
| 472 | 585 | -731.5933519684 | -5.09795 |
| 473 | 585 | -739.7043063279 | -5.09795 |
| 474 | 625 | -817.2443964495 | 16.56113 |
| 475 | 625 | -835.1087734264 | 16.56113 |
| 476 | 625 | -838.4310203320 | 17.20547 |
| 477 | 625 | -831.1137715223 | 17.20547 |
| 478 | 625 | -713.3440268690 | 17.20547 |
| 479 | 625 | -703.5906042517 | 16.56113 |
| 480 | 625 | -710.9078530614 | 16.56113 |
| 481 | 625 | -695.4796498921 | 17.20547 |
| 482 | 729 | -710.7158718413 | 2.4381 |
| 483 | 729 | -723.3892380281 | 2.4381 |
| 484 | 729 | -708.5896338217 | 2.4381 |
| 485 | 729 | -715.7420991849 | 2.4381 |
| 486 | 729 | -729.8894075453 | 2.4381 |
| 487 | 905 | -712.6642748512 | 2.45395 |
| 488 | 905 | -720.7752292108 | 2.45395 |
| 489 | 905 | -729.2884869065 | 2.45395 |
| 490 | 905 | -725.3376410380 | 2.45395 |
| 491 | 905 | -725.9662400009 | 2.45395 |
| 492 | 945 | -777.7788733621 | 9.41812 |
| 493 | 945 | -766.2008433959 | 9.41812 |
| 494 | 945 | -789.0307244532 | 9.41812 |
| 495 | 945 | -778.8742095827 | 9.41812 |
| 496 | 945 | -772.5878630720 | 9.41812 |
| 497 | 945 | -758.0898890364 | 9.41812 |
| 498 | 945 | -660.1574166666 | -41.82726 |
| 499 | 945 | -665.9770264196 | -41.82726 |
| 500 | 945 | -669.2992733252 | -41.82726 |
| 501 | 945 | -665.3484274567 | -41.82726 |
| 502 | 945 | -663.8507884000 | -41.82726 |
| 503 | 945 | -668.4283662324 | -41.82726 |
| 504 | 1625 | -698.3221577328 | -49.36403 |
| 505 | 1625 | -692.1489610632 | -49.36403 |
| 506 | 1625 | -707.6008611156 | -49.36403 |

Table A4.1. 506 invariants for rigid isotopy of 8-crosses.

**Appendix 5**

Let us calculate the antisymmetric matrix $H = P \otimes O$ from Eq.(14) rewriting it with circular permutation of indexes as

$$H_{i,j} H_{i,k} = -[prM(G, \boldsymbol{U})_j]_{k,i} [prM(G, \boldsymbol{U})_k]_{i,j}, \quad (A5.1)$$

or as a linear equation

$$H_{i,j} = -H_{i,k} [prM(G, \boldsymbol{U})_j]_{k,i} [prM(G, \boldsymbol{U})_k]_{i,j} \quad (A5.2)$$

because $|H_{i,k}| = 1$. Now fix $H_{0,1} = 1$ (the other solution is at $H_{0,1} = -1$). For $i = 0$ and $k = 1$ we get from Eq. (A5.2)

$$H_{0,j} = -H_{0,1}[prM(G, \boldsymbol{U})_j]_{1,0}[prM(G, \boldsymbol{U})_1]_{0,j} \equiv -[prM(G, \boldsymbol{U})_j]_{1,0}[prM(G, \boldsymbol{U})_1]_{0,j} \quad (A5.3)$$

that is the $0^{th}$ row of the matrix $H$ and $0^{th}$ column, because $H$ is antisymmetric. Here $j \neq 1$. To calculate the other rows, let us make out of Eq. (A5.2) a recursion relation by choosing $k = i - 1$:

$$H_{i,j} = H_{k,i}[prM(G, \boldsymbol{U})_j]_{k,i}[prM(G, \boldsymbol{U})_k]_{i,j} = H_{i-1,i}[prM(G, \boldsymbol{U})_j]_{i-1,i}[prM(G, \boldsymbol{U})_{i-1}]_{i,j}. \quad (A5.4)$$

Now it is clear that the subsequent rows are calculated from previous ones while taking $i = 1, 2 \ldots, (n-1)$.

**References**

1. P.V. Pikhitsa and S. Pikhitsa, Mutually touching infinite cylinders in the 3d world of lines, Siberian Electronic Mathematical Reports, vol.**16**, 96-120 (2019).
2. P.V. Pikhitsa and S. Pikhitsa, Symmetry, topology and the maximum number of mutually pairwise-touching infinite cylinders: configuration classification, R. Soc. Open Sci., vol. **4**:1, 160729 (2017). (http://dx.doi.org/10.1098/rsos.160729)
3. P.V. Pikhitsa, M. Choi, H.-J. Kim, and S.-H. Ahn, Auxetic lattice of multipods, phys stat solidi b **246**, 2098-2101 (2009). (doi:10.1002/pssb.200982041)
4. O. Ya. Viro and Yu. V. Drobotukhina, Configurations of skew lines, Leningrad Math. J. **1**:4 1027–1050 (1990). arXiv:math/0611374
5. Yu. V. Dobrotukhina, An analogue of the Jones polynomial for links in $RP^3$ and a generalization of the Kauffman-Murasugi theorem, Leningrad Math. J., **2**, N3, 613-630 (1991).
6. V. F. Masurovskii and N.B. Pavlov, Classification of ordered nonsingular configurations of at most seven lines of $RP^3$ up to rigid isotopy, Journal of Mathematical Sciences, vol. **91,** N6, 3508- 3517 (1998).
7. V.F. Mazurovskii, Configuration of six skew lines, Zap. Nauchn. Sem. Leningrad. Otdel. Mat. Inst. Steklov. **167,** 121-134 (1988).
8. V F Mazurovskii , Kauffman polynomials of non-singular configurations projective lines *Russ. Math. Surv.* **44,** 212-213 (1989).
9. A. Connes, *Non commutative geometry*, Academic Press, 1994.

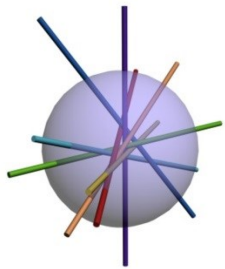